\documentclass[12pt]{article}
\usepackage{amsthm,amsmath,amssymb,amscd,amsfonts,graphicx,bm,xcolor}

\newcommand{\vertiii}[1]{{\left\vert\kern-0.25ex\left\vert\kern-0.25ex\left\vert #1
    \right\vert\kern-0.25ex\right\vert\kern-0.25ex\right\vert}}

\begin{document}

\centerline{\bf Effective Kohn Algorithm for Special Domain}
\centerline{\bf Defined by Functions Depending on All Variables}

\vskip.3in
\bigbreak\centerline{\it Dedicated to the Memory of Jean-Pierre Demailly}

\bigbreak
\centerline{Yum-Tong Siu}

\vskip.3in\bigbreak\noindent{\small\bf Abstract.}  Kohn introduced in 1979 the algorithm of multipliers to study the subelliptc estimate of the $\bar\partial$-Neumann problem for a smooth weakly pseudoconvex domain in a complex Euclidean space which satisfies D'Angelo's finite type condition of a finite bound for the normalized touching order to the boundary for any local possibly singular holomorphic curve in the complex Euclidean space. The problem can be regarded as an example of the formulation of H\"ormander's 1967 hypoelliptic result for the case of complex-valued vector-valued unknowns.  So far the effective solution of Kohn's problem is known only for special domains ${\rm Re}\,z_{n+1}+\sum_{j=1}^N|f_j|^2<0$ in ${\mathbb C}^{n+1}$ with $f_j$ holomorphic in $z_1,\cdots,z_n$, because for such domains it suffices to deal with holomorphic multipliers.  One main obstacle to treat the general smooth case is the need to deal with nonholomorphic multipliers.  This note introduces a new technique to handle nonholomorphic multipliers occurring in more general domains with $f_j$ holomorphic in all the $n+1$ variables $z_1,\cdots,z_n,z_{n+1}$.

\bigbreak\noindent AMS Subject Classification: 32T25, 32T27, 35N15, 13H15.

\bigbreak\noindent Keywords: pseudonconvex domains, finite type domains, multipliers, effective algorithms, multiplicities, Jacobian determinants, fiberwise differentiation, Weierstrass polynomials.

\vskip.3in\bigbreak\noindent{\small\bf\S1.  Introduction.}  In 1967 H\"ormander proved [Ho67, (3.5), Th.4.3 and Th.5.1] that for a linear partial differential operator of the form
$$
P=\sum_{j=1}^r X_j^2+X_0+c
$$
on a domain $\Omega$ in ${\mathbb R}^n$, where $X_0,X_1,\cdots,X_r$ are real smooth vector fields and $c$ is a real-valued smooth function on $\Omega$, if the $m$-fold iterated Lie brackets
$$
X_{j_1},\,[X_{j_1},X_{j_2}],\,  [X_{j_1},\,[X_{j_2},X_{j_3}]],\cdots,\,
[X_{j_1},\,[X_{j_2},\,[X_{j_3},\cdots,X_{j_m}]]]
$$
(for $0\leq j_i\leq r,\,1\leq i\leq m$) span the real tangent space of $\Omega$ at every point, then for a smooth real-valued function $v$ on $\Omega$ with compact support
$$
\|v\|_{\frac{1}{m}}\lesssim\|v\|+\|Pv\|,
$$
where $\|\cdot\|_{\frac{1}{m}}$ is the $\frac{1}{m}$-Sobolev norm on $\Omega$ and $\|\cdot\|$ is the $L^2$ norm on $\Omega$ and the notation $\lesssim$ means that the left-hand side is no more than some constant times the right-hand side.

\bigbreak In 1979 Kohn [Ko79] considered the problem of subellipicity of the $\bar\partial$-Neumann problem on smooth weakly pseudoconvex domains $\Omega$ in ${\mathbb C}^n$.  For a point $P_0$ of the boundary $\partial\Omega$ of $\Omega$, he asked for sufficient conditions for the following subelliptic estimates to hold for some $\varepsilon>0$ and some open neighborhood $U$ of $P_0$ in ${\mathbb C}^n$.
$$
\left\||\varphi|\right\|^2_{\varepsilon}\lesssim\|\bar\partial\varphi\|^2+\|\bar\partial^*\varphi\|^2+\|\varphi\|^2
$$
for any test smooth $(0,q)$-form $\varphi$ with support on $U\cap\bar\Omega$ which is in the domain of the actual adjoint of $\bar\partial$, where $\|\cdot\|$ means the $L^2$ norm over $\Omega$ and $\left\||\cdot|\right\|_{\varepsilon}$ means the Sobolev norm for derivatives of order $\leq\varepsilon$ on $\Omega$ along a smooth family of directions in $\bar\Omega$ which are tangential to the boundary $\partial\Omega$.  The condition for $\varphi$ to be in the domain of the actual adjoint of $\bar\partial$ is that its components with at least one index in the normal direction are zero.  The condition of the vanishing of normal components is the reason for calling the problem the $\bar\partial$-Neumann problem.

\bigbreak In a way, the question of the subelliptic estimates of the $\bar\partial$-Neumann problem can be regarded as an example of the formulation of H\"ormander's 1967 result for the case of complex-valued vector-valued unknowns $\varphi$.  The relation is the closest when $n=2$ and $q=1$, because the complex dimension of the $(1,0)$-tangent space of $\partial\Omega$ is $1$ so that $\varphi$ can be regarded as a complex-valued scalar function.

\bigbreak For Frobenius integrability the formulation in terms of Lie brackets of defining vector fields is equivalent to the dual formulation in terms of exterior differentiation of defining differential forms.  Kohn's approach uses exterior differentiation of differential forms instead of H\"ormander's approach of using Lie brackets of vector fields.  One advantage of Kohn's approach is the root-taking procedure.  As a way of measuring, in H\"ormander's approach, whether the subspace generated by $k$-fold Lie brackets up to a certain $k$ is yet the full tangent space, one can consider the subsheaf of the tangent bundle generated by $k$-fold Lie brackets up to a certain $k$.  In Kohn's approach for the cotangent space, the elements of the corresponding subsheaf are the {\it vector-valued multipliers} $\sigma$ such that its interior product $\sigma\cdot\varphi$ with the test form $\varphi$ satisfies
$$
\left\||\sigma\cdot\varphi|\right\|^2_{\varepsilon}\lesssim\|\bar\partial\varphi\|^2+\|\bar\partial^*\varphi\|^2+\|\varphi\|^2
$$
for some $\varepsilon$ and for all $\varphi$.   When all components of the vector-valued multiplier $\sigma$ are equal, it becomes a {\it scalar multiplier} $f$ such that
$$
\left\||f\varphi|\right\|^2_{\varepsilon}\lesssim\|\bar\partial\varphi\|^2+\|\bar\partial^*\varphi\|^2+\|\varphi\|^2
$$
for some $\varepsilon$ and for all $\varphi$.  From Cramer's rule the determinant of $n-q$ vector-valued multiplier is a scalar multiplier.  Forming the determinant is the same as forming the exterior product.  The root-taking procedure makes it possible to take a root of the absolute-value of a scalar multiplier to get another scalar multiplier for a smaller $\varepsilon$.  More generally, a function dominated by the absolute value of a scalar multiplier is a multiplier.  H\"ormander's procedure of forming the Lie bracket corresponds to Kohn's procedure of taking the $(1,0)$-differential of a scalar multiplier to get a vector-valued multiplier for a smaller $\varepsilon$.

\bigbreak What is new in Kohn's approach is the three procedures of (i) taking the exterior product ({\it i.e.,} the determinant) of vector-valued multipliers to generate a scalar multiplier, (ii) using the domination of a function, in absolute values, by a root of a scalar multiplier to get a scalar multiplier, and (iii) taking the $(1,0)$-exterior differential of a scalar multiplier to generate a vector multiplier.  These procedures used together provide an algorithm of generating scalar multipliers, known as Kohn's algorithm.  Corresponding to H\"ormander's condition of generation of the full tangent space by iterated Lie brackets, Kohn's approach ends up with a subelliptic estimate if the constant function $1$ is generated from Kohn's algorithm as a scalar multiplier.

\bigbreak To avoid the unnecessary notational complicity, we are going to limit ourselves only to the case $q=1$, because the arguments for a general $q$ are rather straightforward modifications of those for the case of $q=1$.  One looks for an easily verifiable geometric condition to guarantee subelliptic estimates for the $\bar\partial$-Neumann problem.
A good condition is the {\it finite type condition} of D'Angelo [DA79] which says that at the prescribed point $P_0$ of $\partial\Omega$ in question, the {\it normalized touching order} of any local holomorphic complex curve $C$ in ${\mathbb C}^n$ to $\partial\Omega$ at $P_0$ is no more than some finite number $p$.  The normalized touching order of $C$ to $\partial\Omega$ at $P_0$ means the supremum of
$$
\frac{\ {\rm ord}_{P_0}(r\circ\psi)\ }{\ \min_{1\leq j\leq n}{\rm ord}_{P_0}\psi_j\ },
$$
when $\Omega=\{r<0\}$ in a neighborhood of $P_0$ in ${\mathbb C}^n$ and $dr\not=0$ at $P_0$, $\psi=(\psi_1,\cdots,\psi_n)$ is a holomorphic map from an open neighborhood of $0$ in ${\mathbb C}$ to ${\mathbb C}^n$ whose image is $C$, and ${\rm ord}_{P_0}$ means the order of vanishing at $P_0$.  The normalized touching order is used, because local complex curves $C$ with singularities at $P_0$ are allowed in the statement of the definition of finite type.

\bigbreak For the case of real-analytic boundary, the work of Diederich-Fornaess [DF78] can be applied to show that D'Angelo's finite type condition implies subelliptic estimate [Ko79, Th.1.19].  However, the result is not effective in the sense that Kohn's algorithm cannot be carried out in $p^*$ steps to yield an $\varepsilon$ as the order of subellipticity such that both $p^*$ and $\varepsilon$ depend explicitly on $p$.  In H\"ormander's setting the number $m$ for the spanning of the tangent space by iterated Lie brackets up to order $m$ plays the role of the finite type $p$ with $\varepsilon=\frac{1}{m}$ as the final order of ellipticity.  One looks for a similar answer in Kohn's situation.  Up to this point the only known effectiveness result is for the following {\it very special domain}
$\Omega$ in ${\mathbb C}^{n+1}$ such that $\tilde U\cap\Omega$ is defined by
$$
r:={\rm Re}(z_{n+1})+\sum_{j=1}^N \left|f\left(z_1,\cdots,z_n\right)\right|^2<0
$$
for some open neighborhood $\tilde U$ of $P_0$ in ${\mathbb C}^{n+1}$, where each of $f_1,\cdots,f_N$ is holomorphic and does not depend on the last variable $z_{n+1}$ (see [Si10, KZ18, KZ21].  This kind of domain is even more special than the {\it special domain} defined in [Ko79, (7.1)] where the holomorphic functions $f_1,\cdots,f_N$ are allowed to depend on the last variable $z_{n+1}$.  For the convenience of referring to the two different kinds of domains, we call the case of defining holomorphic functions independent of $z_{n+1}$ a {\it very special domain} and call the case of defining holomorphic functions depending on $z_{n+1}$ a {\it special domain}.

\bigbreak The distinction between the two different kinds of domains is important, because in the case of independence of $z_{n+1}$, the $(1,0)$ tangent space of $\partial\Omega$ is naturally biholomorphic to the space ${\mathbb C}^n$ with variables $z_1,\cdots,z_n$ so that the arguments about termination can be restricted to scalar multipliers which are holomorphic functions in $z_1,\cdots,z_n$.  For holomorphic scalar multipliers the methods of algebraic geometry can be readily applied.  However, for the case of defining holomorphic functions depending on $z_{n+1}$, one has to consider scalar multipliers which depend on the complex conjugate $\bar z_{n+1}$ of $z_{n+1}$ besides their holomorphic dependence on $z_1,\cdots,z_n, z_{n+1}$, making it impossible to readily apply the methods of algebraic geometry.  In Kohn's 1979 paper, even for the noneffective termination of Kohn's algorithm in [Ko79, Th.7.13] for the case of a special domain with defining holomorphic functions depending on $z_{n+1}$, the theory of Diederich-Fornaess [DF78] needs to be used in the form of [Ko79, Th.6.27] for an intermediate step [Ko79, Prop.7.4] in the argument.

\bigbreak One big obstacle for solving the general case of a smooth weakly pseudoconvex domain of finite type is the method to handle scalar multipliers which depend on the complex conjugate $\bar z_j$ of the complex variables $z_j$ as well as the complex variables $z_j$ themselves.  To provide a way to overcome such an obstacle is the reason why the effective termination of the modified Kohn algorithm for special domains is investigated.  More complicated arguments alow the same lines are still needed for the general smooth.  Here we will not go into the discussion of the general smooth case.

\bigbreak Catlin and D'Angelo presented in [CD10, Prop.4.4] a counter-example of a very special domain in ${\mathbb C}^3$ for which the original Kohn's algorithm of using always the full real radical of the ideal of scalar multipliers at every step cannot be effective. Thus, in order to get effectiveness even for very special domains, we need to modify the original Kohn algorithm to allow the use of an ideal with the property that an effective power of it is contained in the ideal of scalar multipliers in any intermediate step (see [Si10, p.1182, II.1]).  In \S7 we will explain how the noneffectiveness of Kohn's original algorithm in the counter-example of Catlin and D'Angelo becomes effective in the modified Kohn algorithm.  We are going to call the modified Kohn algorithm the {\it effective Kohn algorithm}.

\bigbreak In this note we will prove the effective termination result for the effective Kohn algorithm for the case of special domains where the defining holomorphic functions depend on all the complex variables including $z_{n+1}$.

\bigbreak\noindent{\small\bf\S2. Statement of Result and Topics in Its Argument.}

\bigbreak The statement of our result for the effective termination of the modified Kohn algorithm for special domains is the following.

\bigbreak\noindent(2.1) {\it Main Theorem.}  Let $P_0$ be a point of ${\mathbb C}^{n+1}$ and $f_1,\cdots,f_N$ be holomorphic function germs on ${\mathbb C}^{n+1}$ at $P_0$ which vanish at $0$.
Let $\Omega$ be the domain in ${\mathbb C}^{n+1}$ such that $\tilde U\cap\Omega$ is defined by
$$
r:={\rm Re}(z_{n+1})+\sum_{j=1}^N \left|f\left(z_1,\cdots,z_n,z_{n+1}\right)\right|^2<0
$$
on some open neighborhood $\tilde U$ of $P_0$ in ${\mathbb C}^{n+1}$.  Assume that $\partial\Omega$ is of finite type $\leq p$ at $P_0$ in the sense of D'Angelo.
Then there exists some open neighborhood $U$ of $P_0$ in $\tilde U$ and some positive number $\varepsilon$ such that the subelliptic estimate
$$
\left\||\varphi|\right\|^2_\varepsilon\lesssim\|\bar\partial\varphi\|^2+\|\bar\partial^*\varphi\|^2+\|\varphi\|^2
$$
of order $\varepsilon$ holds for any test smooth $(0,1)$-form $\varphi$ with support on $U\cap\bar\Omega$ which is in the domain of the actual adjoint of $\bar\partial$, where $\|\cdot\|$ means the $L^2$ norm over $\Omega$ and $\left\||\cdot|\right\|_\varepsilon$
means the Sobolev norm for derivatives of order $\leq\varepsilon$ on $\Omega$ along a smooth family of directions in $\bar\Omega$ which are tangential to the boundary $\partial\Omega$.  Moreover, the subelliptic estimate is obtained from an effective Kohn algorithm for multipliers with the steps of generating multipliers depending explicitly and effectively on $p$ so that $\varepsilon$ is explicitly computable from $p$.

\bigbreak\noindent(2.2) {\it Topics in the Argument of the Main Theorem.}  The proof of the Main Theorem involves the following four topics.

\bigbreak\noindent(1) Effective estimation of finite type.

\bigbreak\noindent(2) Multiplicity bound for Jacobian determinant of finite fiber map.

\bigbreak\noindent(3) Algebraic geometric method for termination of effective Kohn algorithm for a finite collection of pre-multipliers.

\bigbreak\noindent(4) New technique to handle non-holomorphic scalar multipliers for special domains where holomorphic defining functions depend on all variables.

\bigbreak\noindent The first three topics were already treated for the case of very special domains.   The reason for treating them here will be explained at the point where each is presented. The fourth topic is the key new technique to prove the Main Theorem.

\bigbreak After the proof of the Main Theorem by the presentation of these four topics, we include only some short remarks on (i) the counter-example of Catlin and D'Angelo on noneffectiveness of the original Kohn algorithm, (ii) the context of the earlier introduction of H\"ormander's iterated Lie bracket conditions by Wei-Liang Chow as a generalization of Carath\'eodory's work on thermodynamics, and (iii) the multitypes of Catlin.

\bigbreak\noindent{\small\bf\S3. Effective Estimation of Finite Type.}  In order to put D'Angelo's finite type condition in a more useful form, we need the following effective description of it given by the following lemma.  The special domain $\Omega$ in Lemma(3.1) below with a prescribed boundary point $P_0$ is as defined in the Main Theorem.

\bigbreak\noindent(3.1) {\it Lemma.} Let $p^*$ be the smallest integer such that
$$
{\mathfrak m}_{{\mathbb C}^{n+1},P_0}^{p^*}\subset\sum_{j=1}^N{\mathcal O}_{{\mathbb C}^{n+1},P_0}f_j,
$$
where ${\mathfrak m}_{{\mathbb C}^{n+1},P_0}$ is the maximum ideal of the local ring ${\mathcal O}_{{\mathbb C}^{n+1},P_0}$ for ${\mathbb C}^{n+1}$ at $P_0$.  Then
$p\leq p^*\leq (n+3)p$.

\bigbreak For the proof of Lemma(3.1) Skoda's result for ideal generation [Sk72, Th.1, pp.555-556] is needed, which is the reason for the factor $n+3$ in the inequality $p\leq p^*\leq (n+3)p$.  This kind of results on ideal generation was treated later in more general settings by Skoda [Sk78, Sk80] and by Demailly [De82].  In our present situation of special domains with $f_j$ depending on $z_{n+1}$, only completely straightforward modifications of the arguments in [Si10, (I.4)] are needed for the proof  which we skip.

\bigbreak\noindent{\small\bf\S4. Multiplicity of Jacobian Determinant for Finite-Fiber Holomorphic Map}.  Let $\pi:{\mathbb C}^{n+1}\to{\mathbb C}^{n+1}$ be the germ of a finite-fiber map at $0$ with $\lambda$ sheets, which maps $0$ to $0$.   More precisely, there exist an open neighborhood $\hat U$ of $0$ in the domain space ${\mathbb C}^{n+1}$ and an open neighborhood $U$ of $0$ in the target space ${\mathbb C}^{n+1}$ such that $\pi$ defines a proper holomorphic map from $\hat U$ to $U$.  We are using ${\mathbb C}^{n+1}$ instead of ${\mathbb C}^n$, because a special domain is defined as a domain in ${\mathbb C}^{n+1}$ and its defining holomorphic functions $f_1,\cdots,f_N$ are holomorphic function germs on ${\mathbb C}^{n+1}$.

\bigbreak\noindent(4.1) {\it Proposition.} The multiplicity of the Jacobian determinant ${\rm Jac}(\pi)$ of the finite-fiber map at $0$ is no more than $\lambda-1$.

\bigbreak\noindent{\it Proof.}  This is the same as [Si10, p.1194, (III.5)] where the standard use of the argument principle in the proof is skipped which is inserted here to make it a little bit easier to follow.  We use $X$ to denote the domain space ${\mathbb C}^{n+1}$ and use $Y$ to denote the target space ${\mathbb C}^{n+1}$.

\bigbreak Let $W$ be an open neighborhood of $0$ in $Y$ and $\tilde W$ be an open neighborhood of $0$ in $X$ such that the map from $\tilde W$ to $W$ induced by $\pi$ is proper with $\lambda$ sheets.  Let $Y_\pi$ be the branching locus in $Y$ for $\pi$.
We can assume without loss of generality (after a linear transformation of the coordinates $z_1,\cdots,z_{n+1}$ in $Y$) that for some open polydisk neighborhood $U=U_1\times\cdots\times U_{n+1}$ of $0$ in $Y$, the set $U_1\times\cdots\times U_n\times\partial U_{n+1}$ is disjoint from $Y_\pi$.  There exists a a subvariety $Y_\pi^\prime$ of in $Y_\pi$ of complex dimension $\leq n-1$ such that for any point $Q$ of $Y_\pi$ not in $Y_\pi^\prime$ and any point $\tilde Q$ in $X$ with $\pi(\tilde Q)$, there exist a local holomorphic coordinate chart $w_1,\cdots,w_{n+1}$ of $Y$ centered at $Q$ and a local holomorphic coordinate chart $\zeta_1,\cdots,\zeta_{n+1}$ of $X$ centered at $\tilde Q$ such that $w_j=\zeta_j$ for $1\leq j\leq n$ and $w_{n+1}=\zeta_{n+1}^{\hat\lambda}$ for some $1\leq\hat\lambda\leq\lambda$.  In other words, at a general point of the branching locus of $\pi$, the map $\pi$ is a power map in the transversal direction.  We can choose a point $P$ of $U_1\times\cdots\times U_n$ such that $(\{P\}\times U_{n+1})$ intersects $Y_\pi$ transversally at $Q_1,\cdots,Q_\ell$ and is disjoint from $Y^\prime_\pi$.  Corresponding to $\hat\lambda$ at $Q$, we have $\hat\lambda_1,\cdots,\hat\lambda_\ell$ at $Q_1,\cdots,Q_\ell$.  Let the disjoint union $\cup_{j=1}^\ell G_j$ be the normalization of the noncompact Riemann surface $\pi^{-1}(\{P\}\times U_{n+1})$.  The boundary $\partial G_j$ of the noncompact Riemann surface $G_j$ is a circle which covers the circle boundary $\{P\}\times\partial U_j$ of $\{P\}\times U_n$ precisely $\lambda_j$ times for $1\leq j\leq\ell$.  The total vanishing order of $z_{n+1}$ on $G_j$ can be computed by applying the argument principle
$$
\frac{1}{2\pi}\Delta_{\partial G_j}{\rm arg}(z_{n+1})=\frac{1}{2\pi i}\int_{G_j}\frac{dz_{n+1}}{z_{n+1}}
$$
to the change of argument $\Delta_{\partial G_j}{\rm arg}(z_{n+1})$ of $z_{n+1}$ along the boundary $\partial G_j$ of $G_j$ which is no more than $\lambda_j$.  This means that the vanishing multiplicity $\sum_{j=1}^\ell(\hat\lambda_j-1)$ of the restriction of ${\rm Jac}(\pi)$ on $(\{P\}\times U_{n+1})$ is no more than $\lambda-1$, because $\lambda=\sum_{j=1}^\ell(\hat\lambda_j)$.

\bigbreak\noindent{\small\bf\S5. Termination of Effective Kohn Algorithm for a Finite Collection of Pre-Multipliers.}  We now treat the third topic, with a new streamlined approach, especially in the argument given in (5.3).
To start out, we need to introduce some terminology.

\bigbreak\noindent(5.1) {\it Terminology.} We start out with a finite collection of holomorphic function germs $f_1,\cdots,f_N$ on ${\mathbb C}^{n+1}$ at a point $P_0$ of ${\mathbb C}^{n+1}$ which without loss of generality we assume to be the origin.  Assume that the ideal generated by $f_1,\cdots,f_N$ contains the $p$-th power of the maximum ideal ${\mathfrak m}_{{\mathbb C}^{n+1},0}$ of ${\mathbb C}^{n+1}$ at the origin for some {\it effective} number $p$.  We call each of $f_1,\cdots,f_N$ a {\it pre-multiplier}.  We introduce the following three kinds of {\it allowable procedures}.

\bigbreak\noindent(i) {\it Jacobian Determinant of Pre-Multipliers.}  For $n+1$ pre-multipliers $g_1,\cdots,g_{n+1}$, their Jacobian determinant
$$
\frac{\partial(g_1,\cdots,g_{n+1})}{\partial(z_1,\cdots,z_{n+1})}
$$
with respect to the coordinates $z_1,\cdots,z_{n+1}$ is a {\it multiplier}.  Every multiplier is a pre-multiplier.  This procedure is allowed only when the germ of the map defined by $g_1,\cdots,g_{n+1}$ from ${\mathbb C}^{n+1}$ to ${\mathbb C}^{n+1}$ at $0$ is a finite-fiber map with an effective number of sheets so that according to Proposition(4.1) the multiplicity of the resulting Jacobian determinant is effective at $0$.  This means that every multiplier so obtained is of effective multiplicity at $0$.  This procedure can also be carried out when each of $g_1,\cdots,g_{n+1}$ is a ${\mathbb C}$-linear combination of $f_1,\cdots,f_N$ under the above condition on the effective number of sheets for the map defined by $g_1,\cdots,g_{n+1}$.  In actual application, appropriate genericity in the choice of ${\mathbb C}$-linear combination is assumed to guarantee the fulfillment of the condition.

\bigbreak\noindent(ii) {\it Elements of Ideal Generated by Finite Number of Multipliers.}  For any finite collection of multipliers $F_1,\cdots,F_m$, any holomorphic function germ $F$ of the form $G_1F_1+\cdots+G_mF_m$ is a multiplier for any holomorphic functions germs $G_1,\cdots,G_m$.  The procedure of constructing a multiplier is allowed only when the new holomorphic function germ $F$ has effective multiplicity at $0$.

\bigbreak\noindent(iii) {\it Replacing Factor of Multiplier by a Root.}  If $f,g$ are holomorphic function germs and $m\in{\mathbb N}$ such that $fg^m$ is a multiplier, then $fg$ is a multiplier.  Since the multiplier $fg^m$ has effective multiplicity at $0$, automatically the number $m$ has to be effective relative to $p$.  This means that the new multiplier $FG$ so obtained is of effective multiplicity at $0$.

\bigbreak\noindent(5.2) {\it Goal of Termination of Effective Kohn Algorithm.}  When we start out with pre-multipliers $f_1,\cdots,f_N$ such that the ideal generated by them contains the $p$-th power of the maximum ideal at the origin with $p$ effective, our goal is to use an effective number of the three kinds of allowable procedures to produce $n+1$ multipliers which generate the maximum ideal of ${\mathbb C}^{n+1}$ at $0$ so that taking their Jacobian would end up with the constant function $1$ as a multiplier.

\bigbreak As the initial step, we select $1\leq i_1<\cdots<i_{n+1}\leq N$ such that the ideal generated by $f_{i_1},\cdots,f_{i_{n+1}}$ contains an effective power of the maximum ideal at $0$.  Let $h_{1,1}$ be the Jacobian determinant of $f_{i_1},\cdots,f_{i_{n+1}}$ with respect to $z_1,\cdots,z_{n+1}$ so that the multiplicity of $h_{1,1}$ is effective at $0$ according to Proposition(4.1).  Choose $h_{1,2},\cdots,h_{1,n}$, each being a generic ${\mathbb C}$-linear combination of $f_1,\cdots,f_N$, such that the ideal generated by $h_{1,1},\cdots,h_{1,n+1}$ contains an effective power of the maximum ideal at $0$.
We now describe an induction process to reach the goal of the termination of effective Kohn algorithm.

\bigbreak\noindent(5.3) {\it Induction Procedure.}   The induction is on $1\leq\nu\leq n+1$ to effectively produce holomorphic function germs $h_{\nu,1},\cdots,h_{\nu,n+1}$ such that (i) $h_{\nu,j}$ is a multiplier for $1\leq j\leq\nu$, (ii) each of $h_{\nu,\nu+1},\cdots,h_{\nu,n+1}$ is a generic ${\mathbb C}$-linear combination of $f_1,\cdots,f_N$, and (iii) the ideal generated by $h_{\nu,1},\cdots,h_{\nu,n+1}$ contains an effective power of the maximum ideal of ${\mathbb C}^{n+1}$ at $0$.  Suppose we have
$h_{\nu,1},\cdots,h_{\nu,n+1}$.  We are going to construct $h_{\nu+1,1},\cdots,h_{\nu+1,n+1}$ by effectively using the three kinds of allowable procedures.
This construction uses fiberwise differentiation and iterates the differential operator defined from the Jacobian determinant, first described in [Si10, (III.7) and (III.8)] for the step of going from $\nu=1$ to $\nu=2$.  The argument of going from the $\nu$-th step to the $(\nu+1)$-step is basically the same, with some obvious modifications.  The application of iterating the differential operator from the Jacobian determinant to provide effectiveness to the example of Catlin and D'Angelo given in \S7 below is a good illustration of this technique.

\bigbreak Since the ideal generated by $h_{\nu,1},\cdots,h_{\nu,n+1}$ contains an effective power of the maximum ideal at $0$, the complex-analytic subspace germ $Y_\nu$ defined by the $\nu$ multipliers $h_{\nu,1},\cdots,h_{\nu,\nu}$ is a complete intersection of complex codimension $\nu$ whose multiplicity at $0$ is effective.
In the arguments that follow, we will use a generic linear coordinates $z_1,\cdots,z_{n+1}$.  To start with, the generic linear coordinates $z_1,\cdots,z_{n+1}$ are chosen such that the projection $\pi_{\nu,j}:(z_1,\cdots,z_{n+1})\mapsto(z_j,z_{\nu+1},\cdots,z_{n+1})$ maps $Y_\nu$ to a complex hypersurface germ $Y_{\nu,j}$ in ${\mathbb C}^{n+2-\nu}$ which is defined by a Weierstrass polynomial $W_{\nu,j}$ in $z_j$ of effective degree $p_\nu$ whose coefficients are functions in the pre-multipliers $h_{\nu,\nu+1},\cdots,h_{\nu,n+1}$ for $1\leq j\leq\nu$.  The vanishing of the Weierstrass polynomial $W_{\nu,j}$ on the complex-analytic subspace $Y_\nu$ implies that $W_{\nu,j}$ belongs to the ideal generated by the $\nu$ multipliers $h_{\nu,1},\cdots,h_{\nu,\nu}$ and is therefore a multiplier with effective multiplicity, according to the second allowable procedure of generating multipliers.

\bigbreak We can write the Jacobian determinant of $W_{\nu,1},\cdots,W_{\nu,\nu},h_{\nu,\nu+1},\cdots,h_{\nu,n+1}$ in the form
$$
\begin{aligned}&(dh_{\nu,\nu+1}\wedge\cdots\wedge dh_{\nu,n+1}\wedge dz_1\wedge\cdots\wedge dz_\nu)\left((\partial_{z_1} W_{\nu,1})\cdots(\partial_{z_\nu}W_{\nu,\nu})\right)\cr
&=\frac{\partial(h_{\nu,\nu+1},\cdots,h_{\nu,n+1})}{\partial(z_{\nu+1},\cdots,z_{n+1})}\left((\partial_{z_1} W_{\nu,1})\cdots(\partial_{z_\nu}W_{\nu,\nu})\right)(dz_1\wedge\cdots\wedge dz_{n+1}).\end{aligned}
$$
It is important to note that for this conclusion the assumption on $W_{\nu,1},\cdots,W_{\nu,\nu}$ can be weakened to the condition that $W_{\nu,j}$ is a polynomial in $z_j,\cdots,z_\nu$ whose coefficients are functions in the pre-multipliers $h_{\nu,\nu+1},\cdots,h_{\nu,n+1}$ for $1\leq j\leq\nu$.
On the $(n+1-\nu)$-dimensional complex-analytic subspace germ defined by $W_{\nu,1},\cdots,W_{\nu,\nu}$,
the function $$\frac{\partial(h_{\nu,\nu+1},\cdots,h_{\nu,n+1})}{\partial(z_{\nu+1},\cdots,z_{n+1})}$$ defines an $(n-\nu)$-dimensional complex-analytic subspace germ which we can project by the map defined by $h_{\nu,\nu+1},\cdots,h_{\nu,n+1}$ to a hypersurface germ in ${\mathbb C}^{n-\nu+1}$ (with coordinates $h_{\nu,\nu+1},\cdots,h_{\nu,n+1}$) defined by $g_\nu(h_{\nu,\nu+1},\cdots,h_{\nu,n+1})$.  As a result,
$$
\begin{aligned}&\left|g_\nu(h_{\nu,\nu+1},\cdots,h_{\nu,n+1})\right|\cr
&\lesssim|W_{\nu,1}|+\cdots+|W_{\nu,\nu}|+\left|\frac{\partial(h_{\nu,\nu+1},\cdots,h_{\nu,n+1})}{\partial(z_{\nu+1},\cdots,z_{n+1})}\right|
\end{aligned}
$$
and
$$
\begin{aligned}&\left|g_\nu(h_{\nu,\nu+1},\cdots,h_{\nu,n})\left((\partial_{z_1} W_{\nu,1})\cdots(\partial_{z_\nu}W_{\nu,\nu})\right)\right|\cr
&\lesssim|W_{\nu,1}|+\cdots+|W_{\nu,\nu}|+\left|\frac{\partial(h_{\nu,\nu+1},\cdots,h_{\nu,n+1})}{\partial(z_{\nu+1},\cdots,z_{n+1})}\left((\partial_{z_1} W_{\nu,1})\cdots(\partial_{z_\nu}W_{\nu,\nu})\right)\right|
\end{aligned}
$$
so that
$$
g_\nu(h_{\nu,\nu+1},\cdots,h_{\nu,n+1})\left((\partial_{z_1} W_{\nu,1})\cdots(\partial_{z_\nu}W_{\nu,\nu})\right)
$$
is a multiplier, because each term on the right-hand side of the above inequality is a multiplier.
This means that the operator
$$
g_\nu(h_{\nu,\nu+1},\cdots,h_{\nu,n+1})((\partial_{z_1} W_{\nu,1})\cdots(\partial_{z_{\nu-1}}W_{\nu,\nu-1}))\partial_{z_\nu}
$$
applied to the polynomial
$W_{\nu,\nu}$ in $z_j,\cdots,z_\nu$ whose coefficients are functions in $h_{\nu,\nu+1},\cdots,h_{\nu,n+1}$ produces a multiplier.  Now we replace
$W_{\nu,\nu}$ by
$$
g_\nu(h_{\nu,\nu+1},\cdots,h_{\nu,n+1})((\partial_{z_1} W_{\nu,1})\cdots(\partial_{z_{\nu-1}}W_{\nu,\nu-1})(\partial_{z_\nu}W_{\nu,\nu}))
$$
to conclude that the operator
$$
(g_\nu(h_{\nu,\nu+1},\cdots,h_{\nu,n+1})(\partial_{z_1} W_{\nu,1})\cdots(\partial_{z_{\nu-1}}W_{\nu,\nu-1}))^2\partial_{z_\nu}^2
$$
applied to the polynomial
$W_{\nu,\nu}$ in $z_j,\cdots,z_\nu$ whose coefficients are functions in $h_{\nu,\nu+1},\cdots,h_{\nu,n+1}$ produces a multiplier.
Again, we apply the same argument, with $W_{\nu,\nu}$ replaced by the polynomial
$$
\left(g_\nu(h_{\nu,\nu+1},\cdots,h_{\nu,n+1})(\partial_{z_1} W_{\nu,1})\cdots(\partial_{z_{\nu-1}}W_{\nu,\nu-1})\right)^2\partial_{z_\nu}^2W_{\nu,\nu}
$$
in $z_j,\cdots,z_\nu$ whose coefficients are functions in $h_{\nu,\nu+1},\cdots,h_{\nu,n+1}$, to produce inductively a multiplier
$$
\left(g_\nu(h_{\nu,\nu+1},\cdots,h_{\nu,n+1})(\partial_{z_1} W_{\nu,1})\cdots(\partial_{z_{\nu-1}}W_{\nu,\nu-1})\right)^{p_\nu}\partial_{z_\nu}^{p_\nu}W_{\nu,\nu}.
$$
Since the degree of $W_{\nu,\nu}$ is a polynomial of degree $p_\nu$ in $z_\nu$, it follows that $\partial_{z_\nu}^{p_\nu}W_{\nu,\nu}$ is a constant.  By taking the $p_\nu$-th root of
$$
\left(g_\nu(h_{\nu,\nu+1},\cdots,h_{\nu,n+1})(\partial_{z_1} W_{\nu,1})\cdots(\partial_{z_{\nu-1}}W_{\nu,\nu-1})\right)^{p_\nu},
$$
we conclude that
$$
g_\nu(h_{\nu,\nu+1},\cdots,h_{\nu,n+1})(\partial_{z_1} W_{\nu,1})\cdots(\partial_{z_{\nu-1}}W_{\nu,\nu-1})
$$
is a multiplier.  Now we apply the same argument to $W_{\nu,\nu-1}$ instead of to $W_{\nu,\nu}$ and use the fact that
the operator
$$
g_\nu(h_{\nu,\nu+1},\cdots,h_{\nu,n})((\partial_{z_1} W_{\nu,1})\cdots(\partial_{z_{\nu-2}}W_{\nu,\nu-2}))\partial_{z_{\nu-1}}
$$
applied to the polynomial
$W_{\nu,\nu-1}$ in $z_j,\cdots,z_{\nu-1}$ whose coefficients are functions in $h_{\nu,\nu+1},\cdots,h_{\nu,n+1}$ produces a multiplier.  Again we replace
$W_{\nu,\nu-1}$ by
$$
g_\nu(h_{\nu,\nu+1},\cdots,h_{\nu,n+1})((\partial_{z_1} W_{\nu,1})\cdots(\partial_{z_{\nu-1}}W_{\nu,\nu-1})(\partial_{z_{\nu-1}}W_{\nu,\nu-1}))
$$
to repeat the same argument to conclude that
$$
g_\nu(h_{\nu,\nu+1},\cdots,h_{\nu,n+1})(\partial_{z_1} W_{\nu,1})\cdots(\partial_{z_{\nu-2}}W_{\nu,\nu-2})
$$
is a multiplier.
We keep on applying the same argument to finally conclude that
$g_\nu(h_{\nu,\nu+1},\cdots,h_{\nu,n+1})$ is a multiplier.  This means that we can set $h_{\nu+1,j}=h_{\nu,j}$ for $1\leq j\leq\nu$ and  $h_{\nu+1,\nu+1}=g_\nu(h_{\nu,\nu+1},\cdots,h_{\nu,n+1})$ and choose each of $h_{\nu+1,\nu+2},\cdots,h_{\nu+1,n+1}$ to be a generic ${\mathbb C}$-linear combination of $f_1,\cdots,f_N$.

\bigbreak At the last step of $\nu=n+1$, each of $h_{n+1,1},\cdots,h_{n+1,n+1}$ is a multiplier and the ideal generated by $h_{n+1,1},\cdots,h_{n+1,n+1}$ contains an effective power of the maximum ideal of ${\mathbb C}^{n+1}$ at $0$.  It means that each $z_1,\cdots,z_{n+1}$ is a multiplier according to the second allowable procedure.  The Jacobian of the $n+1$ multipliers $z_1,\cdots,z_{n+1}$ with respect to the coordinates $z_1,\cdots,z_{n+1}$ is a multiplier which is the constant function $1$.  This effectively terminates the algorithm for a finite collection of pre-multipliers which generate an ideal containing an effective power of the maximum ideal at the origin.

\bigbreak\noindent(5.4){\it Reduction of the Case of Very Special Domain to Effective Termination for Finite Collection of Pre-Multipliers.}  For a very special domain defined by
$$
r={\rm Re}\,z_{n+1}+\sum_{j=1}^N|f_j(z_1,\cdots,z_n)|^2,
$$
with the holomorphic defining functions $f_1,\cdots,f_N$ independent of $z_{n+1}$, one initial scalar multiplier from Kohn's treatment of a smooth weakly pseudoconvex domain [Ko79, Th.1.21] is given by the coefficient of $$(dz_1\wedge\cdots\wedge dz_{n+1})\wedge(d\bar z_1\wedge\cdots\wedge d\bar z_{n+1})$$ in $\partial r\wedge\bar\partial r\wedge(\partial\bar\partial r)^n$, which means that
$$
\sum_{1\leq j_1,\cdots,j_n\leq N}\left|\frac{\partial(f_{j_1},\cdots,f_{j_n})}{\partial(z_1,\cdots,z_n)}\right|^2
$$
is a scalar multiplier.  Each Jacobian determinant
$$
\frac{\partial(f_{j_1},\cdots,f_{j_n})}{\partial(z_1,\cdots,z_n)}
$$
is a holomorphic scalar multiplier, because its absolute-value square is dominated by a scalar multiplier.  Let $F_1,\cdots,F_\ell$ be holomorphic scalar multipliers for $1\leq\ell\leq n$.
Then the coefficient of $$(dz_1\wedge\cdots\wedge dz_{n+1})\wedge(d\bar z_1\wedge\cdots\wedge d\bar z_{n+1})$$ in
$$
\partial r\wedge\bar\partial r\wedge(\partial\bar\partial r)^{n-\ell}\wedge(dF_1\wedge\overline{dF_1}\wedge(dF_\ell\wedge\overline{dF_\ell})
$$
is a scalar multiplier, which means that
$$
\sum_{1\leq j_1,\cdots,j_{n-\ell}\leq N}\left|\frac{\partial(f_{j_1},\cdots,f_{j_{n-\ell}},F_1,\cdots,F_\ell)}{\partial(z_1,\cdots,z_n)}\right|^2
$$
is a scalar multiplier.  Again, each Jacobian determinant
$$
\frac{\partial(f_{j_1},\cdots,f_{j_{n-\ell}},F_1,\cdots,F_\ell)}{\partial(z_1,\cdots,z_n)}
$$
is a holomorphic scalar multiplier, because its absolute-valued square is dominated by a scalar multiplier.  Though we use only some holomorphic scalar multipliers, we can still conclude that the effective Kohn algorithm for a very special domain terminates effectively, because effective procedures applied to them produce the constant function $1$ as a scalar multiplier.

\bigbreak\noindent(5.5){\it Assigned Subellipticity Order for Multipliers in Algebraic Geometric Formulation of Effective Kohn Algorithm.}  The original reason for the study of effective subelliptic estimates is to come up with some explicit description of the order of subellipticity at the end when the constant function $1$ is finally shown to be a scalar multiplier.  For that purpose we would like to remark on the order of subellipticity assigned to a scalar multiplier obtained at each step.

\bigbreak For the smooth case of a weakly pseudoconvex domain defined by $\{r<0\}$ with $dr\not=0$, the orders of subellipticity for multipliers are given in Kohn's 1979 paper as follows.

\bigbreak\noindent(i) The order of subellipticity for the scalar multiplier $r$ is $1$ (see [Ko79, p.94, (4.8)]).

\bigbreak\noindent(ii) The order of subellipticity is $\frac{1}{2}$ for the vector-multiplier which is the $(1,0)$-form on the boundary of the domain obtained by evaluating the Levi-form $\partial\bar\partial r$ on any $(0,1)$-vector field tangential to the boundary of the domain (see [Ko79, p.96, (4.17)]).

\bigbreak\noindent(iii) The order of subellipticity for the scalar multiplier $f$ is $\frac{\delta}{m}$ if the order of subellipticity of the scalar multiplier $|f|^m$ is $\delta$ (see [Ko79, p.98, (4.35)]).

\bigbreak\noindent(iv) If the order of subellipticity for the scalar multiplier $f$ is $\delta$, then the order of subellipticity for the vector multiplier $\partial f$ is $\frac{\delta}{2}$ {see [Ko79, p.99, (4.42)]).

\bigbreak Translated into holomorphic multipliers constructed from holomorphic pre-multipliers by allowable procedures, the assigned orders of subellipticity are as follows.

\bigbreak\noindent(i) The assigned order of subellipticity of the Jacobian determinant formed from some of the initial functions $f_1,\cdots,f_N$ and some multipliers $F_1,\cdots,F_\ell$ generated in the algorithm is $\frac{1}{2}$ times the minimum of $1$ and the assigned order of subellipticity of $F_1,\cdots,F_\ell$.

\bigbreak\noindent(ii) The assigned order of subellipticity for $FG$ is $\frac{1}{m}$ times the assigned order of subellipticity for $F^mG$.

\bigbreak\noindent{\small\bf\S6. New Techniques to Treat the Special Domain with Defining Holomorphic Functions Depend on All Variables}.  For a special domain with defining holomorphic functions depend on all the variables, when we use Kohn's original construction of scalar multipliers, we end up with scalar multipliers which are not holomorphic.  The new techniques introduced to treat such a special domain seek to reduce the problem to the effective termination of the algorithm for a finite collection of holomorphic pre-multipliers, which we already solved as the third topic.  For this reduction we use orthogonal local frames of $(1,0)$-forms adapted
to the boundary $\partial\Omega$ in some of our formulations.

\bigbreak\noindent(6.1) {\it Levi Form and Use of Exterior Product with $\partial r\wedge\bar\partial r$.}  Let $\omega_1,\cdots,\omega_n,\omega_{n+1}$ be an orthogonal local frames of $(1,0)$-forms such that the pullback of $\omega_{n+1}$ to $\partial\Omega$ is identically zero.  This means that $\omega_{n+1}$ is parallel to $\partial r$.  Write $\partial\bar\partial r=\sum_{i,j=1}^{n+1} c_{i\bar j}\,\omega_i\wedge\overline{\omega_j}$ so that the Levi-form is given by $\left(c_{i\bar j}\right)_{1\leq i,k\leq n-1}$.  An alternative description is to use orthonormal vector fields $L_1,\cdots,L_n,L_{n+1}$ with $L_1,\cdots,L_n$ along the tangential direction and $T=L_{n+1}-\overline{L_{n+1}}$ so that we have $c_{i\bar j}=\langle\partial\bar\partial r,\,L_i\wedge\overline{L_j}\rangle$ or
$$
[L_i,\overline{L_j}]=c_{i\bar j}T+\sum_{k=1}^na_{i\bar j}^kL_k+\sum_{k=1}^nb_{i\bar j}^{\bar k}\overline{L_k},
$$
where $a_{i\bar j}$ and $b_{i\bar j}$ are scalar functions.
A smooth test $(0,1)$-form $\varphi$ is in the domain of actual adjoint $\bar\partial^*$ if and only if the normal component of $\varphi$ vanishes, which means $\varphi=\sum_{i=1}^n\varphi_{\bar i}\,\overline{\omega_k}$.  The key property is [Ko79, Proposition(4.7)(C)], which states that the $(1,0)$-form which is the interior product ${\rm int}(\theta)\partial\bar\partial r$ of the Levi-form $\partial\bar\partial r$ with any smooth $(0,1)$-form $\theta$ with $\left<\theta,\partial r\right>\equiv 0$ is a vector multiplier.  The important point is that the condition $\left<\theta,\partial r\right>\equiv 0$ needs to hold in order for ${\rm int}(\theta)\partial\bar\partial r$ to be a vector multiplier.

\bigbreak Suppose $F_1,\cdots,F_\ell$ are multipliers.  Let $\partial F_j=\sum_{k=1}^{n+1} F_{jk}\omega_k$.  Then
$${\mathcal D}_{\mu_1,\cdots,\mu_{n-\ell}}:=\det\left(\begin{matrix}c_{1,\overline{\mu_1}}&c_{2,\overline{\mu_1}}&\cdots&c_{n,\overline{\mu_1}}\cr
c_{1,\overline{\mu_2}}&c_{2,\overline{\mu_2}}&\cdots&c_{n,\overline{\mu_2}}\cr
\cdot&\cdot&\cdots&\cdot\cr
\cdot&\cdot&\cdots&\cdot\cr
\cdot&\cdot&\cdots&\cdot\cr
c_{1,\overline{\mu_{n-\ell}}}&c_{2,\overline{\mu_{n-\ell}}}&\cdots&c_{n,\overline{\mu_{n-\ell}}}\cr
F_{1,1}&F_{1,2}&\cdots&F_{1,n}\cr
F_{2,1}&F_{2,2}&\cdots&F_{2,n}\cr\cdot&\cdot&\cdots&\cdot\cr
\cdot&\cdot&\cdots&\cdot\cr
\cdot&\cdot&\cdots&\cdot\cr
F_{\ell,1}&F_{\ell,2}&\cdots&F_{\ell,n}\cr
\end{matrix}\right)$$
is a multiplier for any $1\leq\mu_1<\cdots<\mu_{n-\ell}\leq n$, which means that the coefficient ${\mathcal C}$ of $\bigwedge_{j=1}^{n+1}(dz_j\wedge d\bar z_j)$ in
$$
(\partial r\wedge\bar\partial r)\wedge(\partial\bar\partial r)^{n-\ell}\wedge\left(\bigwedge_{j=1}^\ell(\partial F_j\wedge\bar\partial F_j)\right)
$$
is a multiplier, because
$$
|{\mathcal C}|\lesssim\sum_{1\leq\mu_1<\cdots<\mu_{n-\ell}\leq n}\left|{\mathcal D}_{\mu_1,\cdots,\mu_{n-\ell}}\right|^2.
$$
The order of subellipticity for each ${\mathcal D}_{\mu_1,\cdots,\mu_{n-\ell}}$ is $\frac{1}{2}$ times the minimum of $1$ and the order of subellipticity of each $F_1,\cdots,F_\ell$.

\bigbreak\noindent(6.2){\it Computation of Some Scalar Multipliers for the Special Domain.}  The special domain which we consider is
$$
r(z,\bar z)=z_{n+1}+\overline{z_{n+1}}+\sum_{j=1}^N f_j(z_1,\cdots,z_n,z_{n+1})\,\overline{f_j(z_1,\cdots,z_n,z_{n+1})}.
$$
The boundary point is some $P_0$ in $\partial\Omega$.  We are going to use the condition that $f_j(z_1,\cdots,z_n,z_{n+1})$ is a holomorphic function germ on ${\mathbb C}^{n+1}$ at $P_0$ which vanishes at $P_0$.  Note that we cannot impose that condition that $P_0$ is the origin, because $r$ has to vanish at $P_0$ and for notational simplicity we do not want to add the constant $-\sum_{j=1}^N|f_j(0)|^2$ to $r$ as in [Ko79, (7.1)].
We have
$$
\begin{aligned}
\partial r&=dz_{n+1}+\sum_{\mu=1}^N(df_\mu(z))\overline{f_\mu(z)},\cr
\bar\partial r&=d\bar z_{n+1}+\sum_{\mu=1}^N f_\mu(z)\overline{(df_\mu(z))},\cr
\end{aligned}
$$
so that
$$
\begin{aligned}
\partial r\wedge\bar\partial r&=\left(dz_{n+1}+\sum_{\mu=1}^N(df_\mu(z))\overline{f_\mu(z)}\right)\wedge\left(d\bar z_{n+1}+\sum_{\mu=1}^N f_\mu(z)\overline{(df_\mu(z))}\right)\cr
&=dz_{n+1}\wedge d\bar z_{n+1}+\sum_{\mu=1}^N (f_\mu(z)\,dz_{n+1})\wedge\overline{(df_\mu(z))}\cr
&+\sum_{\mu=1}^N\overline{(f_\mu(z)\,dz_{n+1})}\wedge(df_\mu(z))
+\sum_{\mu,\nu=1}^N(f_\nu(z)\,df_\mu(z))\wedge\overline{(f_\mu(z)\,df_\nu(z))}.\cr
\end{aligned}
$$
We also have
$$
\partial\bar\partial r=\sum_{\mu=1}^N df_\mu(z)\wedge \overline{df_\mu(z)}
$$
so that
$$
\begin{aligned}\frac{1}{n!}(-1)^{\frac{n(n-1)}{2}}(\partial\bar\partial r)^n&=\frac{1}{n!}(-1)^{\frac{n(n-1)}{2}}\left(\sum_{\mu=1}^N df_\mu(z)\wedge \overline{df_\mu(z)}\right)^n\cr
&=\sum_{1\leq\mu_1<\cdots<\mu_n\leq N}\left(df_{\mu_1}\wedge\cdots\wedge df_{\mu_n}(z)\right)\wedge\overline{\left(df_{\mu_1}\wedge\cdots\wedge df_{\mu_n}(z)\right)}.\cr
\end{aligned}
$$
From
{\small$$
\begin{aligned}
&\frac{1}{n!}(-1)^{\frac{n(n-1)}{2}}\partial r\wedge\bar\partial r\wedge(\partial\bar\partial r)^n\cr
&=
\bigg(dz_{n+1}\wedge d\bar z_{n+1}+\sum_{\mu=1}^N (f_\mu(z)\,dz_{n+1})\wedge\overline{(df_\mu(z))}\cr
&+\sum_{\mu=1}^N\overline{(f_\mu(z)\,dz_{n+1})}\wedge(df_\mu(z))
+\sum_{\mu,\nu=1}^N(f_\nu(z)\,df_\mu(z))\wedge\overline{(f_\mu(z)\,df_\nu(z))}\bigg)\wedge\cr
&\bigg(\sum_{1\leq\mu_1<\cdots<\mu_n\leq N}\left(df_{\mu_1}(z)\wedge\cdots\wedge df_{\mu_n}(z)\right)\wedge\overline{\left(df_{\mu_1}\wedge\cdots\wedge df_{\mu_n}(z)\right)}\bigg)\cr
&=
dz_{n+1}\wedge d\bar z_{n+1}\wedge
\sum_{1\leq\mu_1<\cdots<\mu_n\leq N}\left(df_{\mu_1}(z)\wedge\cdots\wedge df_{\mu_n}(z)\right)\wedge\overline{\left(df_{\mu_1}\wedge\cdots\wedge df_{\mu_n}(z)\right)}\cr
&
+\sum_{\mu=1}^N(f_\mu(z)\,dz_{n+1})\wedge\overline{(df_\mu(z))}\wedge
\sum_{1\leq\mu_1<\cdots<\mu_n\leq N}\left(df_{\mu_1}(z)\wedge\cdots\wedge df_{\mu_n}(z)\right)\wedge\overline{\left(df_{\mu_1}\wedge\cdots\wedge df_{\mu_n}(z)\right)}
\cr
&+\sum_{\mu=1}^N\overline{(f_\mu(z)\,dz_{n+1})}\wedge(df_\mu(z))\wedge\sum_{1\leq\mu_1<\cdots<\mu_n\leq N}\left(df_{\mu_1}\wedge\cdots\wedge df_{\mu_n}(z)\right)\wedge\overline{\left(df_{\mu_1}\wedge\cdots\wedge df_{\mu_n}(z)\right)}\cr
&+\sum_{\mu,\nu=1}^N(f_\nu(z)df_\mu(z))\wedge\overline{(df_\nu(z)f_\mu(z))}\wedge
\sum_{1\leq\mu_1<\cdots<\mu_n\leq N}\left(df_{\mu_1}(z)\wedge\cdots\wedge df_{\mu_n}(z)\right)\wedge\overline{\left(df_{\mu_1}(z)\wedge\cdots\wedge df_{\mu_n}(z)\right)}\cr
\end{aligned}
$$
}

\noindent it follows that the coefficient of $$dz_1\wedge\cdots\wedge dz_{n+1}\wedge d\bar z_1\wedge\cdots\wedge d\bar z_{n+1}$$ in
$$\partial r\wedge\bar\partial r\wedge(\partial\bar\partial r)^n$$ (up to a sign) is equal to
$$
\begin{aligned}
&
\sum_{1\leq\mu_1<\cdots<\mu_n\leq N}\left(\frac{\partial(f_{\mu_1},\cdots,f_{\mu_n})}{\partial(z_1,\cdots,z_n)}\right)
\overline{\left(\frac{\partial(f_{\mu_1},\cdots,f_{\mu_n})}{\partial(z_1,\cdots,z_n)}\right)}\cr
&
+\sum_{\mu=1}^N\sum_{1\leq\mu_1<\cdots<\mu_n\leq N}\left(\frac{\partial(f_{\mu_1},\cdots,f_{\mu_n})}{\partial(z_1,\cdots,z_n)}\,f_\mu\right)
\overline{\left(\frac{\partial(f_{\mu_1},\cdots,f_{\mu_n},f_\mu)}{\partial(z_1,\cdots,z_{n+1})}\right)}\cr
&+\sum_{\mu=1}^N\sum_{1\leq\mu_1<\cdots<\mu_n\leq N}\left(\frac{\partial(f_{\mu_1},\cdots,f_{\mu_n},f_\mu)}{\partial(z_1,\cdots,z_{n+1})}\right)
\overline{\left(\frac{\partial(f_{\mu_1},\cdots,f_{\mu_n})}{\partial(z_1,\cdots,z_n)}f_\mu\right)}\cr
&+\sum_{\mu,\nu=1}^N\sum_{1\leq\mu_1<\cdots<\mu_n\leq N}\left(\frac{\partial(f_{\mu_1},\cdots,f_{\mu_n},f_\mu)}{\partial(z_1,\cdots,z_{n+1})}f_\nu\right)
\overline{\left(\frac{\partial(f_{\mu_1},\cdots,f_{\mu_n},f_\nu)}{\partial(z_1,\cdots,z_{n+1})}f_\mu\right)},\cr
\end{aligned}
$$
which can be rewritten as
$$
\sum_{1\leq\mu_1<\cdots<\mu_n\leq N}
\left|\frac{\partial(f_{\mu_1},\cdots,f_{\mu_n})}{\partial(z_1,\cdots,z_n)}+\sum_{\mu=1}^N
\frac{\partial(f_{\mu_1},\cdots,f_{\mu_n},f_\mu)}{\partial(z_1,\cdots,z_{n+1})}\overline{f_\mu}\right|^2.
$$
This is a smooth multiplier so that each summand is a smooth multiplier.  Thus, each
$$
\frac{\partial(f_{\mu_1},\cdots,f_{\mu_n})}{\partial(z_1,\cdots,z_n)}+\sum_{\mu=1}^N
\frac{\partial(f_{\mu_1},\cdots,f_{\mu_n},f_\mu)}{\partial(z_1,\cdots,z_{n+1})}\overline{f_\mu}
$$
is a smooth multiplier.  The order of subellipticity of
$$
\frac{\partial(f_{\mu_1},\cdots,f_{\mu_n})}{\partial(z_1,\cdots,z_n)}+\sum_{\mu=1}^N
\frac{\partial(f_{\mu_1},\cdots,f_{\mu_n},f_\mu)}{\partial(z_1,\cdots,z_{n+1})}\overline{f_\mu}
$$
is $\frac{1}{2}$.

\bigbreak Choose an open ball neighborhood $U_\delta$ of $P_0$, which is centered at $P_0$ with radius $\delta>0$ sufficiently small, and let $M$ be the maximum of the absolute value of the first order derivatives of $f_\mu$ for $1\leq\mu\leq N$ and let $\eta$ be the maximum of the absolute value of $f_\mu$ for $1\leq\mu\leq N$.  Note that the vanishing of $f_\mu$ at $P_0$ means that $\eta$ can be made smaller than any prescribed positive number by choosing $\delta>0$ sufficiently small.  Here we use the notation that $\widehat{f_{\mu_j}}$ means the removal of $f_{\mu_j}$. Similar notations will also be used later.  By the Laplace expansion of
$$
\frac{\partial(f_{\mu_1},\cdots,f_{\mu_n},f_\mu)}{\partial(z_1,\cdots,z_{n+1})}
$$
of the determinant of order $n+1$ with respect to the last column consisting of $\frac{\partial f_\nu}{\partial z_{n+1}}$ for $\nu=\mu_1,\cdots,\mu_n,\mu$ as entries,
we can write
$$
\begin{aligned}
\frac{\partial(f_{\mu_1},\cdots,f_{\mu_n},f_\mu)}{\partial(z_1,\cdots,z_{n+1})}
&=-\frac{\partial(f_{\mu_1},\cdots,f_{\mu_n})}{\partial(z_1,\cdots,z_n)}\,\frac{\partial f_\mu}{\partial z_{n+1}}\cr
&\quad+\sum_{j=1}^n(-1)^{n-j}\frac{\partial(f_{\mu_1},\cdots,\widehat{f_{\mu_j}},\cdots,f_{\mu_n},f_\mu)}{\partial(z_1,\cdots,z_n)}\,\frac{\partial f_{\mu_j}}{\partial z_{n+1}},
\end{aligned}
$$
from which we obtain the estimate
$$
\left|\frac{\partial(f_{\mu_1},\cdots,f_{\mu_n},f_\mu)}{\partial(z_1,\cdots,z_{n+1})}\right|
\leq M\bigg(\left|\frac{\partial(f_{\mu_1},\cdots,f_{\mu_n})}{\partial(z_1,\cdots,z_n)}\right|
+\sum_{j=1}^n\left|\frac{\partial(f_{\mu_1},\cdots,\widehat{f_{\mu_j}},\cdots,f_{\mu_n},f_\mu)}{\partial(z_1,\cdots,z_n)}\right|\bigg)
$$
and
$$
\left|\frac{\partial(f_{\mu_1},\cdots,f_{\mu_n},f_\mu)}{\partial(z_1,\cdots,z_{n+1})}\overline{f_\mu}\right|
\leq M\eta\bigg(\left|\frac{\partial(f_{\mu_1},\cdots,f_{\mu_n})}{\partial(z_1,\cdots,z_n)}\right|
+\sum_{j=1}^n\left|\frac{\partial(f_{\mu_1},\cdots,\widehat{f_{\mu_j}},\cdots,f_{\mu_n},f_\mu)}{\partial(z_1,\cdots,z_n)}\right|\bigg)
$$
when $z$ is in the open neighborhood $U_\delta$ of the origin so that
$$
\begin{aligned}&\sum_{1\leq\mu_1<\cdots<\mu_n\leq N}\left|\frac{\partial(f_{\mu_1},\cdots,f_{\mu_n})}{\partial(z_1,\cdots,z_n)}\right|\cr
&\leq
\sum_{1\leq\mu_1<\cdots<\mu_n\leq N}\left|\frac{\partial(f_{\mu_1},\cdots,f_{\mu_n})}{\partial(z_1,\cdots,z_n)}+\sum_{\mu=1}^N
\frac{\partial(f_{\mu_1},\cdots,f_{\mu_n},f_\mu)}{\partial(z_1,\cdots,z_{n+1})}\overline{f_\mu}\right|\cr
&\quad+
\sum_{1\leq\mu_1<\cdots<\mu_n\leq N}\left|\sum_{\mu=1}^N
\frac{\partial(f_{\mu_1},\cdots,f_{\mu_n},f_\mu)}{\partial(z_1,\cdots,z_{n+1})}\overline{f_\mu}\right|\cr
&\leq
\sum_{1\leq\mu_1<\cdots<\mu_n\leq N}\left|\frac{\partial(f_{\mu_1},\cdots,f_{\mu_n})}{\partial(z_1,\cdots,z_n)}+\sum_{\mu=1}^N
\frac{\partial(f_{\mu_1},\cdots,f_{\mu_n},f_\mu)}{\partial(z_1,\cdots,z_{n+1})}\overline{f_\mu}\right|\cr
&\quad+
\sum_{1\leq\mu_1<\cdots<\mu_n\leq N}M\eta\bigg(\left|\frac{\partial(f_{\mu_1},\cdots,f_{\mu_n})}{\partial(z_1,\cdots,z_n)}\right|
+\sum_{j=1}^n\left|\frac{\partial(f_{\mu_1},\cdots,\widehat{f_{\mu_j}},\cdots,f_{\mu_n},f_\mu)}{\partial(z_1,\cdots,z_n)}\right|\bigg)\cr
&\leq
\sum_{1\leq\mu_1<\cdots<\mu_n\leq N}\left|\frac{\partial(f_{\mu_1},\cdots,f_{\mu_n})}{\partial(z_1,\cdots,z_n)}+\sum_{\mu=1}^N
\frac{\partial(f_{\mu_1},\cdots,f_{\mu_n},f_\mu)}{\partial(z_1,\cdots,z_{n+1})}\overline{f_\mu}\right|\cr
&\quad+
(n+1)M\eta\sum_{1\leq\mu_1<\cdots<\mu_n\leq N}\left|\frac{\partial(f_{\mu_1},\cdots,f_{\mu_n})}{\partial(z_1,\cdots,z_n)}\right|\cr
\end{aligned}
$$
and
$$
\begin{aligned}&\sum_{1\leq\mu_1<\cdots<\mu_n\leq N}\left|\frac{\partial(f_{\mu_1},\cdots,f_{\mu_n})}{\partial(z_1,\cdots,z_n)}\right|\cr
&\leq
\frac{1}{1-(n+1)M\eta}\sum_{1\leq\mu_1<\cdots<\mu_n\leq N}\left|\frac{\partial(f_{\mu_1},\cdots,f_{\mu_n})}{\partial(z_1,\cdots,z_n)}+\sum_{\mu=1}^N
\frac{\partial(f_{\mu_1},\cdots,f_{\mu_n},f_\mu)}{\partial(z_1,\cdots,z_{n+1})}\overline{f_\mu}\right|,\cr
\end{aligned}
$$
when $\delta>0$ is chosen to satisfy $(n+1)M\eta<1$.
This means that each
$$
\frac{\partial(f_{\mu_1},\cdots,f_{\mu_n})}{\partial(z_1,\cdots,z_n)}
$$
is a holomorphic multiplier at $P_0$.  Its order of subellipticity is $\frac{1}{2}$.
From
$$
\left|\frac{\partial(f_{\mu_1},\cdots,f_{\mu_n},f_\mu)}{\partial(z_1,\cdots,z_{n+1})}\right|
\leq M\bigg(\left|\frac{\partial(f_{\mu_1},\cdots,f_{\mu_n})}{\partial(z_1,\cdots,z_n)}\right|
+\sum_{j=1}^n\left|\frac{\partial(f_{\mu_1},\cdots,\widehat{f_{\mu_j}},\cdots,f_{\mu_n},f_\mu)}{\partial(z_1,\cdots,z_n)}\right|\bigg)
$$
for $z$ in the open neighborhood $U_\delta$ of $P_0$, it follows that each
$$
\frac{\partial(f_{\mu_1},\cdots,f_{\mu_n},f_\mu)}{\partial(z_1,\cdots,z_{n+1})}
$$
is also a holomorphic multiplier at $P_0$.  Its order of subellipticity is $\frac{1}{2}$.  This paves the way for us to use the algebraic geometric method in the third topic for the effective termination for a finite collection of holomorphic pre-multipliers.

\bigbreak\noindent(6.3){\it Computation of More Scalar Multipliers by Induction Argument.}  Fix $\ell\in{\mathbb N}\cup\{0\}$.  Assume that we have newly generated holomorphic multipliers $F_1,\cdots,F_\ell$.  We would like to show that by induction on $\ell$ that
$$
\frac{\partial(f_{\mu_1},\cdots,f_{\mu_{n-\ell}},F_1,\cdots,F_\ell)}{\partial(z_1,\cdots,z_n)}
$$
is a holomorphic multiplier at $P_0$.  We start with the induction hypothesis that the statement is true for $\ell$ and we would like to verify the statement when $\ell$ is replaced by $\ell+1$ so long as $\ell+1\leq n$.

\bigbreak Let $F$ be a newly generated multiplier
The argument below does not care how $F$ arises as a holomorphic multiplier, as long as $dF$ is a holomorphic vector-multiplier.  We use the notation $\pm$ to avoid the need to write down an explicit power of $-1$.
From
$$
\begin{aligned}
&\pm\partial r\wedge\bar\partial r\wedge(\partial\bar\partial r)^{n-\ell-1}\wedge\left(\bigwedge_{\lambda=1}^\ell(dF_\lambda\wedge\overline{dF_\lambda})\right)\wedge(dF\wedge\overline{dF})\cr
&=
\bigg(dz_{n+1}\wedge d\bar z_{n+1}+\sum_{\mu=1}^N (f_\mu(z)\,dz_{n+1})\wedge\overline{(df_\mu(z))}\cr
&+\sum_{\mu=1}^N\overline{(f_\mu(z)\,dz_{n+1})}\wedge(df_\mu(z))
+\sum_{\mu,\nu=1}^N(f_\nu(z)\,df_\mu(z))\wedge\overline{(f_\mu(z)\,df_\nu(z))}\bigg)\wedge\cr
&\bigg(\sum_{1\leq\mu_1<\cdots<\mu_{n-\ell-1}\leq N}\left(df_{\mu_1}\wedge\cdots\wedge df_{\mu_{n-\ell-1}}\wedge dF_1\wedge\cdots\wedge dF_\ell\wedge dF\right)\cr
&\qquad\qquad\qquad\wedge\overline{\left(df_{\mu_1}\wedge\cdots\wedge df_{\mu_{n-\ell-1}}\wedge  dF_1\wedge\cdots\wedge dF_\ell\wedge dF\right)}\bigg)\cr
&=
dz_{n+1}\wedge d\bar z_{n+1}\wedge
\sum_{1\leq\mu_1<\cdots<\mu_{n-\ell-1}\leq N}\left(df_{\mu_1}\wedge\cdots\wedge df_{\mu_{n-\ell-1}}\wedge dF_1\wedge\cdots\wedge dF_\ell\wedge dF\right)\cr
&\qquad\qquad\qquad\wedge\overline{\left(df_{\mu_1}\wedge\cdots\wedge df_{\mu_{n-\ell-1}}\wedge dF_1\wedge\cdots\wedge dF_\ell\wedge dF\right)}\cr
&
+\sum_{\mu=1}^N(f_\mu\,dz_{n+1})\wedge\overline{(df_\mu)}
\wedge\sum_{1\leq\mu_1<\cdots<\mu_{n-\ell-1}\leq N}\left(df_{\mu_1}\wedge\cdots\wedge df_{\mu_{n-\ell-1}}\wedge dF_1\wedge\cdots\wedge dF_\ell\wedge dF \right)\cr
&\qquad\qquad\qquad\wedge\overline{\left(df_{\mu_1}\wedge\cdots\wedge df_{\mu_{n-\ell-1}}\wedge dF_1\wedge\cdots\wedge dF_\ell\wedge dF\right)}
\cr
&+\sum_{\mu=1}^N\overline{(f_\mu\,dz_{n+1})}\wedge(df_\mu)\wedge\sum_{1\leq\mu_1<\cdots<\mu_{n-\ell-1}\leq N}\left(df_{\mu_1}\wedge\cdots\wedge df_{\mu_{n-\ell-1}}\wedge dF_1\wedge\cdots\wedge dF_\ell\wedge dF\right)\cr
&\qquad\qquad\qquad\wedge\overline{\left(df_{\mu_1}\wedge\cdots\wedge df_{\mu_{n-\ell-1}}\wedge dF_1\wedge\cdots\wedge dF_\ell\wedge dF\right)}\cr
&+\sum_{\mu,\nu=1}^N(f_\nu\, df_\mu)\wedge\overline{(df_\nu f_\mu)}\wedge
\sum_{1\leq\mu_1<\cdots<\mu_{n-\ell-1}\leq N}\left(df_{\mu_1}\wedge\cdots\wedge df_{\mu_{n-\ell-1}}\wedge dF_1\wedge\cdots\wedge dF_\ell\wedge dF\right)\cr
&\qquad\qquad\qquad\wedge\overline{\left(df_{\mu_1}\wedge\cdots\wedge df_{\mu_{n-\ell-1}}\wedge dF_1\wedge\cdots\wedge dF_\ell\wedge dF\right)}\cr
\end{aligned}
$$
it follows that the coefficient of $$dz_1\wedge\cdots\wedge dz_{n+1}\wedge d\bar z_1\wedge\cdots\wedge d\bar z_{n+1}$$ in
$$\partial r\wedge\bar\partial r\wedge(\partial\bar\partial r)^{n-\ell}\wedge\left(\bigwedge_{\lambda=1}^\ell(dF_\lambda\wedge\overline{dF_\lambda})\right)\wedge(dF\wedge\overline{dF})$$ (up to a sign) is equal to
{\small$$
\begin{aligned}
&
\sum_{1\leq\mu_1<\cdots<\mu_{n-\ell-1}\leq N}\left(\frac{\partial(f_{\mu_1},\cdots,f_{\mu_{n-\ell-1}},F_1,\cdots,F_\ell,F)}{\partial(z_1,\cdots,z_n)}\right)
\overline{\left(\frac{\partial(f_{\mu_1},\cdots,f_{\mu_{n-\ell-1}},F_1,\cdots,F_\ell,F)}{\partial(z_1,\cdots,z_n)}\right)}\cr
&
+\sum_{\mu=1}^N\sum_{1\leq\mu_1<\cdots<\mu_{n-\ell-1}\leq N}\left(\frac{\partial(f_{\mu_1},\cdots,f_{\mu_{n-\ell-1}},F_1,\cdots,F_\ell,F)}{\partial(z_1,\cdots,z_n)}\,f_\mu\right)
\overline{\left(\frac{\partial(f_{\mu_1},\cdots,f_{\mu_{n-\ell-1}},F_1,\cdots,F_\ell,F)}{\partial(z_1,\cdots,z_{n+1})}\right)}\cr
&+\sum_{\mu=1}^N\sum_{1\leq\mu_1<\cdots<\mu_{n-\ell-1}\leq N}\left(\frac{\partial(f_{\mu_1},\cdots,f_{\mu_{n-\ell-1}},F_1,\cdots,F_\ell,F,f_\mu)}{\partial(z_1,\cdots,z_{n+1})}\right)
\overline{\left(\frac{\partial(f_{\mu_1},\cdots,f_{\mu_{n-\ell-1}},F)}{\partial(z_1,\cdots,z_n)}f_\mu\right)}\cr
&+\sum_{\mu,\nu=1}^N\sum_{1\leq\mu_1<\cdots<\mu_{n-\ell-1}\leq N}\left(\frac{\partial(f_{\mu_1},\cdots,f_{\mu_{n-\ell-1}},F_1,\cdots,F_\ell,F)}{\partial(z_1,\cdots,z_{n+1})}f_\nu\right)
\overline{\left(\frac{\partial(f_{\mu_1},\cdots,f_{\mu_{n-\ell-1}},F)}{\partial(z_1,\cdots,z_{n+1})}f_\mu\right)},\cr
\end{aligned}
$$
}

\noindent which can be rewritten as
{\small$$
\sum_{1\leq\mu_1<\cdots<\mu_{n-\ell-1}\leq N}
\left|\frac{\partial(f_{\mu_1},\cdots,f_{\mu_{n-\ell-1}},F_1,\cdots,F_\ell,F)}{\partial(z_1,\cdots,z_n)}+\sum_{\mu=1}^N
\frac{\partial(f_{\mu_1},\cdots,f_{\mu_{n-\ell-1}},F_1,\cdots,F_\ell,F,f_\mu)}{\partial(z_1,\cdots,z_{n+1})}\overline{f_\mu}\right|^2.
$$
}

\noindent This is a smooth multiplier so that each summand is a smooth multiplier.  Thus, each
$$
\frac{\partial(f_{\mu_1},\cdots,f_{\mu_{n-\ell-1}},F_1,\cdots,F_\ell,F)}{\partial(z_1,\cdots,z_n)}+\sum_{\mu=1}^N
\frac{\partial(f_{\mu_1},\cdots,f_{\mu_{n-\ell-1}},F_1,\cdots,F_\ell,F,f_\mu)}{\partial(z_1,\cdots,z_{n+1})}\overline{f_\mu}
$$
is a smooth multiplier.
Choose an open ball neighborhood $U_\delta$ of $P_0$, which is centered at $P_0$ with radius $\delta>0$ sufficiently small. Let $M$ be the maximum of the absolute value of the first order derivatives of $F$ and $f_\mu$ for $1\leq\mu\leq N$ and let $\eta$ be the maximum of the absolute value of $f_\mu$ for $1\leq\mu\leq N$.  By the Laplace expansion of
$$
\frac{\partial(f_{\mu_1},\cdots,f_{\mu_{n-\ell-1}},F_1,\cdots,F_\ell,F,f_\mu)}{\partial(z_1,\cdots,z_{n+1})}
$$
of the determinant of order $n+1$ with respect to the last column consisting of $\frac{\partial f_\nu}{\partial z_{n+1}}$ (for $\nu=\mu_1,\cdots,\mu_{n-1},\mu$) and $\frac{\partial F_j}{\partial z_{n+1}}$ (for $1\leq j\leq\ell$) and $\frac{\partial F}{\partial z_{n+1}}$ as entries,
we can write
$$
\begin{aligned}
&\frac{\partial(f_{\mu_1},\cdots,f_{\mu_{n-\ell-1}},F_1,\cdots,F_\ell,F,f_\mu)}{\partial(z_1,\cdots,z_{n+1})}\cr
&=-\frac{\partial(f_{\mu_1},\cdots,f_{\mu_{n-\ell-1}},F_1,\cdots,F_\ell,F)}{\partial(z_1,\cdots,z_n)}\,\frac{\partial f_\mu}{\partial z_{n+1}}\cr
&\quad+\sum_{j=1}^{n-\ell-1}(-1)^{n-j}\frac{\partial(f_{\mu_1},\cdots,\widehat{f_{\mu_j}},\cdots,f_{\mu_{n-\ell-1}},F_1,\cdots,F_\ell,F,f_\mu)}{\partial(z_1,\cdots,z_n)}\,\frac{\partial f_{\mu_j}}{\partial z_{n+1}}\cr
&\quad+\sum_{j=1}^\ell(-1)^{\ell+j-1}\frac{\partial(f_{\mu_1},\cdots,f_{\mu_{n-\ell-1}},F_1,\cdots,\widehat{F_j},\cdots,F_\ell,f_\mu)}{\partial(z_1,\cdots,z_n)}\,\frac{\partial F_j}{\partial z_{n+1}},
\end{aligned}
$$
from which we obtain the estimate
$$
\begin{aligned}&\left|\frac{\partial(f_{\mu_1},\cdots,f_{\mu_{n-\ell-1}},F_1,\cdots,F_\ell,F,f_\mu)}{\partial(z_1,\cdots,z_{n+1})}\right|\cr
&\leq M\bigg(\left|\frac{\partial(f_{\mu_1},\cdots,f_{\mu_{n-\ell-1}},F_1,\cdots,F_\ell,F)}{\partial(z_1,\cdots,z_n)}\right|\cr
&\qquad+
\sum_{j=1}^{n-\ell-1}\left|\frac{\partial(f_{\mu_1},\cdots,\widehat{f_{\mu_j}},\cdots,f_{\mu_{n-\ell-1}},F_1,\cdots,F_\ell,F,f_\mu)}{\partial(z_1,\cdots,z_n)}\right|\cr
&\qquad+\sum_{j=1}^\ell\left|\frac{\partial(f_{\mu_1},\cdots,f_{\mu_{n-\ell-1}},F_1,\cdots,\widehat{F_j},\cdots,F_\ell,f_\mu)}{\partial(z_1,\cdots,z_n)}\right|\bigg)
\end{aligned}
$$
and
$$
\begin{aligned}&\left|\frac{\partial(f_{\mu_1},\cdots,f_{\mu_{n-\ell-1}},F_1,\cdots,F_\ell,F,f_\mu)}{\partial(z_1,\cdots,z_{n+1})}\overline{f_\mu}\right|\cr
&\leq M\eta\bigg(\left|\frac{\partial(f_{\mu_1},\cdots,f_{\mu_{n-\ell-1}},F_1,\cdots,F_\ell,F)}{\partial(z_1,\cdots,z_n)}\right|\cr
&\qquad+
\sum_{j=1}^{n-\ell-1}\left|\frac{\partial(f_{\mu_1},\cdots,\widehat{f_{\mu_j}},\cdots,f_{\mu_{n-\ell-1}},F_1,\cdots,F_\ell,F,f_\mu)}{\partial(z_1,\cdots,z_n)}\right|\cr
&\qquad+\sum_{j=1}^\ell\left|\frac{\partial(f_{\mu_1},\cdots,f_{\mu_{n-\ell-1}},F_1,\cdots,\widehat{F_j},\cdots,F_\ell,f_\mu)}{\partial(z_1,\cdots,z_n)}\right|\bigg)\cr
\end{aligned}
$$
when $z$ is in the open neighborhood $U_\delta$ of $P_0$ so that
{\small$$
\begin{aligned}&\sum_{1\leq\mu_1<\cdots<\mu_{n-\ell-1}\leq N}\left|\frac{\partial(f_{\mu_1},\cdots,f_{\mu_{n-\ell-1}},F_1,\cdots,F_\ell,F)}{\partial(z_1,\cdots,z_n)}\right|\cr
&\leq
\sum_{1\leq\mu_1<\cdots<\mu_{n-\ell-1}\leq N}\left|\frac{\partial(f_{\mu_1},\cdots,f_{\mu_{n-\ell-1}},F_1,\cdots,F_\ell,F)}{\partial(z_1,\cdots,z_n)}+\sum_{\mu=1}^N
\frac{\partial(f_{\mu_1},\cdots,f_{\mu_{n-\ell-1}},F_1,\cdots,F_\ell,F,f_\mu)}{\partial(z_1,\cdots,z_{n+1})}\overline{f_\mu}\right|\cr
&\quad+
\sum_{1\leq\mu_1<\cdots<\mu_{n-\ell-1}\leq N}\left|\sum_{\mu=1}^N
\frac{\partial(f_{\mu_1},\cdots,f_{\mu_{n-\ell-1}},F_1,\cdots,F_\ell,F,f_\mu)}{\partial(z_1,\cdots,z_{n+1})}\overline{f_\mu}\right|\cr
&\leq
\sum_{1\leq\mu_1<\cdots<\mu_{n-\ell-1}\leq N}\left|\frac{\partial(f_{\mu_1},\cdots,f_{\mu_{n-\ell-1}},F_1,\cdots,F_\ell,F)}{\partial(z_1,\cdots,z_n)}+\sum_{\mu=1}^N
\frac{\partial(f_{\mu_1},\cdots,f_{\mu_{n-\ell-1}},F_1,\cdots,F_\ell,F,f_\mu)}{\partial(z_1,\cdots,z_{n+1})}\overline{f_\mu}\right|\cr
&\quad+
\sum_{1\leq\mu_1<\cdots<\mu_{n-\ell-1}\leq N}M\eta\bigg(\left|\frac{\partial(f_{\mu_1},\cdots,f_{\mu_{n-\ell-1}},F_1,\cdots,F_\ell,F)}{\partial(z_1,\cdots,z_n)}\right|\cr
&\qquad+
\sum_{j=1}^{n-\ell-1}\left|\frac{\partial(f_{\mu_1},\cdots,\widehat{f_{\mu_j}},\cdots,f_{\mu_{n-\ell-1}},F_1,\cdots,F_\ell,F,f_\mu)}{\partial(z_1,\cdots,z_n)}\right|\cr
&\qquad+\sum_{j=1}^\ell\left|\frac{\partial(f_{\mu_1},\cdots,f_{\mu_{n-\ell-1}},F_1,\cdots,\widehat{F_j},\cdots,F_\ell,f_\mu)}{\partial(z_1,\cdots,z_n)}\right|\bigg)\cr
&\leq
\sum_{1\leq\mu_1<\cdots<\mu_{n-\ell-1}\leq N}\left|\frac{\partial(f_{\mu_1},\cdots,f_{\mu_{n-\ell-1}},F_1,\cdots,F_\ell,F)}{\partial(z_1,\cdots,z_n)}+\sum_{\mu=1}^N
\frac{\partial(f_{\mu_1},\cdots,f_{\mu_{n-\ell-1}},F_1,\cdots,F_\ell,F,f_\mu)}{\partial(z_1,\cdots,z_{n+1})}\overline{f_\mu}\right|\cr
&\qquad+
(n-\ell+1)M\eta\sum_{1\leq\mu_1<\cdots<\mu_{n-\ell-1}\leq N}\left|\frac{\partial(f_{\mu_1},\cdots,f_{\mu_{n-\ell-1}},F_1,\cdots,F_\ell,F)}{\partial(z_1,\cdots,z_n)}\right|\cr
&\qquad+\sum_{1\leq\mu_1<\cdots<\mu_{n-\ell-1}\leq N}M\eta\,
\sum_{j=1}^\ell\left|\frac{\partial(f_{\mu_1},\cdots,f_{\mu_{n-\ell-1}},F_1,\cdots,\widehat{F_j},\cdots,F_\ell,f_\mu)}{\partial(z_1,\cdots,z_n)}\right|\cr
\cr
\end{aligned}
$$
}

\noindent and
$$
\begin{aligned}&\sum_{1\leq\mu_1<\cdots<\mu_{n-\ell-1}\leq N}\left|\frac{\partial(f_{\mu_1},\cdots,f_{\mu_{n-\ell-1}},F_1,\cdots,F_\ell,F)}{\partial(z_1,\cdots,z_n)}\right|\cr
&\leq
\frac{1}{1-(n-\ell+1)M\eta}\sum_{1\leq\mu_1<\cdots<\mu_{n-\ell-1}\leq N}\bigg|\frac{\partial(f_{\mu_1},\cdots,f_{\mu_{n-\ell-1}},F_1,\cdots,F_\ell,F)}{\partial(z_1,\cdots,z_n)}\cr
&\qquad\qquad\qquad\quad+\sum_{\mu=1}^N
\frac{\partial(f_{\mu_1},\cdots,f_{\mu_{n-\ell-1}},F_1,\cdots,F_\ell,F,f_\mu)}{\partial(z_1,\cdots,z_{n+1})}\overline{f_\mu}\bigg|\cr
&\qquad+\frac{M\eta}{1-(n-\ell+1)M\eta}\sum_{1\leq\mu_1<\cdots<\mu_{n-\ell-1}\leq N}\sum_{j=1}^\ell\left|\frac{\partial(f_{\mu_1},\cdots,f_{\mu_{n-\ell-1}},F_1,\cdots,\widehat{F_j},\cdots,F_\ell,f_\mu)}{\partial(z_1,\cdots,z_n)}\right|\cr
\end{aligned}
$$
when $\delta>0$ is chosen to satisfy $(n-\ell+1)M\eta<1$.
This means that each
$$
\frac{\partial(f_{\mu_1},\cdots,f_{\mu_{n-\ell-1}},F_1,\cdots,F_\ell,F)}{\partial(z_1,\cdots,z_n)}
$$
is a holomorphic multiplier at $P_0$.  From
$$
\begin{aligned}&\left|\frac{\partial(f_{\mu_1},\cdots,f_{\mu_{n-\ell-1}},F_1,\cdots,F_\ell,F,f_\mu)}{\partial(z_1,\cdots,z_{n+1})}\right|\cr
&\leq M\bigg(\left|\frac{\partial(f_{\mu_1},\cdots,f_{\mu_{n-\ell-1}},F_1,\cdots,F_\ell,F)}{\partial(z_1,\cdots,z_n)}\right|\cr
&\qquad\qquad+\left|\frac{\partial(f_{\mu_1},\cdots,f_{\mu_{n-\ell-1}},f_\mu)}{\partial(z_1,\cdots,z_n)}\right|
+\sum_{j=1}^{n-1}\left|\frac{\partial(f_{\mu_1},\cdots,\widehat{f_{\mu_j}},\cdots,f_{\mu_{n-\ell-1}},F_1,\cdots,F_\ell,F,f_\mu)}{\partial(z_1,\cdots,z_n)}\right|\bigg)
\end{aligned}
$$
for $z$ in the open neighborhood $U_\delta$ of $P_0$, it follows that each
$$
\frac{\partial(f_{\mu_1},\cdots,f_{\mu_{n-\ell-1}},F_1,\cdots,F_\ell,F,f_\mu)}{\partial(z_1,\cdots,z_{n+1})}
$$
is also a holomorphic multiplier at $P_0$.  This finishes the induction argument.  The order of subellipticity of
$$
\frac{\partial(f_{\mu_1},\cdots,f_{\mu_{n-\ell-1}},F_1,\cdots,F_\ell,F,f_\mu)}{\partial(z_1,\cdots,z_{n+1})}
$$
is $\frac{1}{2}$ times the minimum of $1$ and the order of subellipticity of $F_1,\cdots,F_\ell,F$.

\bigbreak\noindent(6.4){\it Effective Termination of Kohn Algorithm for Special Domain.}  By Lemma(3.2) we know that the holomorphic function germs $f_1,\cdots,f_N$ on ${\mathbb C}^{n+1}$ at $P_0$ can be used as a finite collection of pre-multipliers on ${\mathbb C}^{n+1}$ at $P_0$.  From (6.2) we know that the Jacobian determinant
$$
\frac{\partial(f_{\mu_1},\cdots,f_{\mu_n},f_\mu)}{\partial(z_1,\cdots,z_{n+1})}
$$
is a holomorphic multiplier at $P_0$.  Moreover, for newly generated holomorphic multipliers $F_1,\cdots,F_\ell$, the induction argument in (6.3) shows that
$$
\frac{\partial(f_{\mu_1},\cdots,f_{\mu_{n-\ell+1}},F_1,\cdots,F_\ell)}{\partial(z_1,\cdots,z_{n+1})}
$$
is again a holomorphic multiplier at $P_0$.  We can now use the effective termination of the algorithm for the finite collection of holomorphic multipliers $f_1,\cdots,f_N$ on ${\mathbb C}^{n+1}$ at $P_0$ to conclude that effectively we can produce the constant function $1$ as a scalar multiplier for the special domain $\Omega$ at the boundary point $P_0$.

\bigbreak\noindent{\small\bf\S7. Remark on Example of Catlin and D'Angelo.}  Catlin and D'Angelo [CD10,p.81, Prop.4.4] considered the very special domain in ${\mathbb C}^3$ defined by the two pre-multipliers $F_1(z,w)=z^M$ and $F_2(z,w)=w^N+wz^K$ with $K>M\geq 2$, where $M$ and $N$ are effective but $K$ is not.  For this example, Kohn's original algorithm of taking the full real radical at every step is not effective.  However, the modified Kohn algorithm of using fiberwise differentiation and iterating the differential operator defined from the Jacobian determinant is effective.  The purpose of this remark is to use this counter-example to illustrate the key point of the modified Kohn algorithm on iterating the differential operator from the Jacobian determinant, as presented in (5.3).

\bigbreak The reason for the noneffectiveness of Kohn's algorithm for this example is as follows.  Since there are two initial pre-multipliers $F_1$ and $F_2$ in two complex variables $z,w$, the first step is to form the Jacobian determinant
$$
h=\det\left(\begin{matrix}Mz^{M-1}&0\cr Kwz^{K-1}&Nw^{N-1}+z^K\cr\end{matrix}\right)
=Mz^{M-1}(Nw^{N-1}+z^K).$$
Its full radical is the principal ideal $I_0$ generated by $g=z(Nw^{N-1}+z^K)$.  The next step is to add to the principal ideal the two Jacobian determinants
$$
\frac{\partial(F_1,g)}{\partial(z,w)}\quad{\rm and}\quad\frac{\partial(F_2,g)}{\partial(z,w)}
$$
to form the ideal $J_1$ and then its full radical $I_1$.  Since
$$
\frac{\partial(F_2,g)}{\partial(z,w)}
$$
contains the factor $z^K$, the function $z^K$ belongs to $J_1$.  However, $z^{K-1}$ does not belong to $J_1$, because the three generators of $J_1$ modulo $w=0$ all contain the factor $z^K$.  Thus, to go from $J_1$ to $I_1$ the $K$-th root of the element $z^K$ of $J_1$ has to be taken and $K$ is not effective.

\bigbreak Let us see the reason for the effectiveness of our modified Kohn algorithm of using fiberwise differentiation and iterating the differential operator defined from the Jacobian determinant.  We consider $$h=\frac{\partial(F_1,F_2)}{\partial(z,w)}$$ as the result of applying the differential operator
$$
L=-(\partial_w F_1)\partial_z+(\partial_z F_1)\partial_w=Mz^{M-1}\partial_w
$$
to the function $F_2$.  We can iterate the application of $L$ to get multipliers
$$h_k=L^kF_2=M^kN(N-1)\cdots(N-k+1)z^{k(M-1)}w^{N-k}$$
for $k\geq 2$.  By taking the the effective $N(M-1)$-th root of $h_N$, we conclude that $z$ is a multiplier and $wz^K$ is a multiplier.  Then $$w^N=(w^N+wz^K)-wz^K=F_2-wz^K$$ is a pre-multiplier.  The Jacobian determinant
$$
\det\left(\begin{matrix}1&0\cr 0&Nw^{N-1}\cr\end{matrix}\right)=Nw^{N-1}
$$
for $z$ and $w^N$ is a multiplier.
By taking its effective $(N-1)$-th root of $w^{N-1}$, we conclude that $w$ is a multiplier. Finally $1$ which is the Jacobian determinant of the multipliers $z,w$ with respect to the coordinates $z,w$ is a multiplier.  Thus the modified Kohn algorithm effectively terminates.

\bigbreak\noindent{\small\bf\S8. Caratheodory's Axiomatic Approach to Entropy in Thermodynamics and the Iterated Lie Bracket Condition.}  Since H\"ormander's hypoellipticiy result on sums of squares of vector fields was mentioned at the beginning of this note as an important starting point for the investigation of subelliptic estimates, we would like to make some historical remarks.  The iterated Lie bracket condition used in H\"ormander hypoelliptic results on sums of squares of vector fields was first by Wei-Leung Chow in 1939 [Ch39] in another related context.  Chow's work generalized Caratheodory's result [Ca09, p.369] on the integrability condition for the differential $1$-form
$$
\omega:=dx_0+X_1dx_1+X_2dx_2+\cdots+X_ndx_n
$$
in ${\mathbb R}^{n+1}$ with coordinates $x_0,x_1,\cdots,x_n$.  Carath\'eodory proved that $\omega$ is automatically integrable when the {\it nonreachability condition} is satisfied that in each neighborhood of an arbitrary point $P$ there exists some point which is not reachable from $P$ along broken curves on which the pullback of the $1$-form $\omega$ vanishes.  Caratheodory's argument is as follows and uses mainly geometry rather than analysis.

\bigbreak Let $P_i$ be a sequence of points approaching $P$ which are not reachable from $P$. Let $G$ be the line through $P$ parallel to the $x_0$-axis.  First, we would like to argue that without loss of generality we can choose a subsequence of $P_i$ to be on $G$.  Suppose the contrary so that we can assume each $P_i$ to be outside $G$. Let $\Pi_i$ be the $2$-plane containing $G$ and $P_i$ and let $C_i$ be the integral curve through $P_i$ for the pullback of $\omega$ to $\Pi_i$.  Let $Q_i$ be the intersection of $G$ and $C_i$ and we know that $Q_i$ (being reachable from $P_i$) must be different from $P$ but approaches $P$ as $P_i$ approaches $P$.  This means that $P_i$ can be chosen to be on $G$ after we replace $P_i$ by $Q_i$.

\bigbreak Take another line $G_1$ which is different from $G$ and is parallel to the $x_0$-axis. Let $S$ be an arbitrary $2$-dimensional cylinder with axis parallel to the $x_0$-axis such that $G$ and $G_1$ are inside $S$.  We pullback $\omega$ to $S$ and consider the integral curve through $P$ and let $M_{S,P}$ be the point of its intersection with $G_1$.  We claim that $M_{S,P}$ is independent of $S$, otherwise as $S$ moves to $\hat S$ we have $M_{\hat S, Q}=M_{S,P}$ for some $Q$ on $G$ close to $P$ which can reach $P$ through the two integral curves containing the point $M_{\hat S, Q}=M_{S,P}$.  Moreover, $\hat S$ can be chosen to make $Q$ cover an open neighborhood of $P$ in $G$, which is a contradiction.

\bigbreak We now write $M_{G_1,P}$ instead of $M_{S,P}$, because the dependence on $S$ is only through the condition that $S$ contains $G_1$.  We can vary $G_1$ with $n$ parameters to get $M_{G_1,P}$ varying with $n$ parameters.  This $n$-fold containing all $M_{G_1,P}$ as $G_1$ varies is the integral $n$-fold of $\omega$ which contains $P$.

\bigbreak Coupled with Frobenius result, Caratheodory's result can be reformulated as the equivalence of the condition of nonreachability and the condition that the Lie brackets of the vector fields in the kernel of $\omega$ do not add new vector field not in the kernel of $\omega$.

\bigbreak Carath\'eodory's work sought to axiomatization of the theory of entropy.  In the equation for $\omega$ the term $dx_0$ is the infinitesimal change of work done on a thermal system.  For example,  $X_i$ can mean the pressure of a gas and $dx_i$ can mean the infinitesmal change of the volume of a gas.  A point $P$ can mean the state of a thermal system.  The unreachability condition means that for a given state $P$ there are states $Q$ nearby not reachable from $P$.  Then automatically the $1$-form $\omega$ is integrable to yield integrable hypersurfaces.  All states on the same integrable hypersurface have the same ``entropy'' and can be reachable from one another. His conclusion is that the notion of ``entropy'' arises from the axiomatic assumption of the existence of nearby unreachable states.

\bigbreak In the setting of Carath\'eodory, a metric between two points $P$ and $Q$ in ${\mathbb R}^{n+1}$ can be defined by using only paths which are broken curves on which the pullback of $1$-form $\omega$ vanishes.  Such a metric is called a Carnot-Carath\'eodory metric, because of Carnot's work [Ca24] on the Carnot cycle for the most efficient heat engine under the laws of thermodynamics, especially the reversibility from conservation of entropy.

\bigbreak\noindent{\small\bf\S9. Chow's Generalization Caratheodory's Result.}  Consider a system $B_r$ defined by a subbundle of the tangent space of constant rank $r$.  Suppose iterated Lie brackets extend $B_r$ to $B_s$ with $s=r+m$.  Then $m$ is called the {\it index} for the system $B_r$.  In any neighborhood of any point $P$, the set of reachable points by broken curves with tangent in the subbundle $B_r$ forms a submanifold of real dimension $s=r+m$.  In Caratheodory's original result, $m=0$ and $r=n$ in the ambient space of real dimension $n+1$ with coordinates $x_0,x_1,\cdots,x_n$ so that in any neighborhood of a point the submanifold of reachable points is of dimension $n=r+m$
and cannot contain the whole neighborhood.  Chow's arguments are mainly analytical instead of geometrical.  He related the condition of the spanning of the full tangent space by iterated Lie brackets to the condition of reachability.  In particular, he used the Lie bracket $[\xi,\eta]$ as the infinitesimal generator for the group $e^{\xi t}e^{\eta t}e^{-\xi t}e^{-\eta t}$, which essentially is how H\"ormander used the iterated Lie brackets, except that H\"ormander used the broken integral curves between two points to compute the H\"older norm from the difference quotient with a denominator of fractional power as well as a smoothing process.

\bigbreak Rashevski in [Ra1938] independently obtained a result similar to Chow's.

\bigbreak\noindent{\small\bf\S10. Catlin's Approach of Multitype.}  We would like to point out that to study the problem of subelliptic estimates for the $\bar\partial$-Neumann problem for smooth weakly pseudoconvex domains, Catlin [Ca83, Ca84, Ca87] introduced the notion of multitype which carries more information than D'Angelo's finite type.
We will not go into Catlin's approach here.

\bigbreak\noindent{\bf References}

\medbreak\noindent[Ca24]  Sadi Carnot, {\it R\'eflexions sur la puissance motrice du feu et sur les machines propres \`a d\'evelopper cette puissance}. Paris: Chez Bachelier Libraire, 1824.

\medbreak\noindent[Ca09] Constantin Carath\'eodory,
Untersuchungen \"uber die Grundlagen der Thermodynamik. {\it Math. Ann.} \textbf{67} (1909), 355 -- 386.

\medbreak\noindent[Ca83] D. Catlin,
Necessary conditions for subellipticity of the $\bar\partial$-Neumann problem. {\it Ann. of Math.} {\bf 117} (1983),
147 -- 171.

\medbreak\noindent[Ca84] D. Catlin,
Boundary invariants of
pseudoconvex domains. {\it Ann. of Math.} {\bf 120} (1984),
529 -- 586.

\medbreak\noindent[Ca87] D. Catlin,
Subelliptic estimates for the $\bar\partial$-Neumann problem on pseudoconvex domains. {\it Ann. of Math.} {\bf 126} (1984),
131 -- 191.

\medbreak\noindent[CD10]
David W. Catlin and John P. D'Angelo,
Subelliptic estimates.  {\it Complex analysis}, pp. 75 -- 94,
Trends Math., Birkh\"auser/Springer Basel AG, Basel, 2010.

\medbreak\noindent[Ch39]
W. L. Chow, \"Uber Systeme von linearen partiellen Differentialgleichungen erster Ordnung, {\it Math. Ann.}
{\bf 117} (1939), 98 -- 105.

\medbreak\noindent[DA79] J. P. D'Angelo,
Finite type conditions for
real hypersurfaces. {\it J. Diff. Geom.} {\bf 14} (1979), 59 -- 66.

\medbreak\noindent[DA82] J. P. D'Angelo,
Real hypersurfaces, orders of contact, and applications. {\it Ann. of Math.} \textbf{115} (1982), 615 -- 637.

\medbreak\noindent[DA92] J. P. D'Angelo,
Finite type conditions and subelliptic estimates. In: {\it Modern methods in complex analysis} (Princeton, NJ, 1992), volume \textbf{137} of Ann. of Math. Stud., pages 63 -- 78.

\medbreak\noindent[DA93] J. P. D'Angelo,
Several complex variables and the geometry of real hypersurfaces. {\it Studies in Advanced Mathematics}. CRC Press, Boca Raton, FL, 1993.

\medbreak\noindent[De82]
Jean-Pierre Demailly, Estimations pour l'op\'erateur d'un fibr\'e vectoriel holomorphe semi-positif au-dessus d'une vari\'et\'e k\"ahl\'erienne compl\`ete.
{\it Ann Scient \'Ecole Norm. Sup.} \textbf{15} (1982), 457 -- 511.  p.481, Th.6.2

\medbreak\noindent[DF78] K. Diederich and J. E. Fornaess,
Pseudoconvex domains with real-analytic boundary. {\it Ann. of
Math.} {\bf 107} (1978), 371 -- 384.

\medbreak\noindent[Ho67]] Lars H\"ormander, Lars
Hypoelliptic second order differential equations. {\it Acta Math.} \textbf{119} (1967), 147 -- 171.

\medbreak\noindent[KZ18] Sung-Yeon Kim and Dmitri Zaitsev. Jet vanishing orders and effectivity of Kohn's algorithm in dimension 3. {\it Asian J. Math.} \textbf{22} (2018), 545 -- 568.

\medbreak\noindent[KZ21] Sung-Yeon Kim and Dmitri Zaitsev. Triangular resolutions and effectiveness for holomorphic subelliptic multipliers. {\it Adv. Math.} \textbf{387}, (2021).

\medbreak\noindent[Ko79] J.~J.~Kohn, Subellipticity of the $\bar
\partial $-Neumann problem on pseudo-convex domains: sufficient
conditions. {\it Acta Math.} {\bf 142} (1979), 79 -- 122.

\medbreak\noindent[Ra38]
P. K. Rashevsky, Any two points of a totally nonholonomic space may be connected by an admissible
line, {\it Uch. Zap. Ped. Inst. im. Liebknechta, Ser. Phys. Math.} {\bf 2} (1938), 83 -- 94.

\medbreak\noindent[Si10]
Yum-Tong Siu. Effective termination of Kohn's algorithm for subelliptic multipliers. {\it Pure Appl. Math. Q.} \textbf{6} (2010), Special Issue: In honor of Joseph J. Kohn. Part 2: 1169--1241.

\medbreak\noindent [Sk72] H.  Skoda, Application des techniques
$L^2$ \`a la th\'eorie des id\'eaux d'une alg\`ebre de fonctions
holomorphes avec poids, {\it Ann. Sci. Ec. Norm. Sup.} {\bf 5}
(1972), 548 -- 580.

\medbreak\noindent [Sk78] H.  Skoda,
Morphismes surjectifs de fibr\'es vectoriels semi-positifs.  {\it Ann. Sci. Ec. Norm. Sup.} {\bf 11}
(1978), 577 -- 611.

\medbreak\noindent [Sk80] H.  Skoda, Rel\`evement des sections globales dans les fibres semi-positifs. S\'eminaire P. Lelong-H. Skoda (Analyse),
19\textsuperscript{e} ann\'ee, 1978-1979, {\it Lecture Notes in Math.} \textbf{822}, Springer-Verlag, Berlin, Heidelberg,
New York, 1980.

\bigbreak\noindent{\it Author's mailing address:} Department of
Mathematics, Harvard University, Cambridge, MA 02138, U.S.A.

\medbreak\noindent{\it Author's e-mail address:}
siu@math.harvard.edu
\end{document}